\documentclass[11pt,reqno]{amsart}
\usepackage{enumerate}
\usepackage{amsfonts}
\usepackage{amsmath}
\usepackage{amsthm}
\usepackage{amssymb}
\usepackage{latexsym}
\usepackage{multicol}
\usepackage{verbatim}
\usepackage{amscd}
\usepackage{hyperref}
\usepackage{amsmath, amscd}
\usepackage{pb-diagram}

\newtheorem*{cor}{Corollary}
\newtheorem*{lem}{Lemma}
\newtheorem*{prop}{Proposition}

\theoremstyle{definition}
\newtheorem*{defn}{Definition}
\theoremstyle{definition}
\newtheorem*{thm}{Theorem}

\newenvironment{pf}{\proof}{\endproof}
\newcounter{cnt}
\newenvironment{enumerit}{\begin{list}{{\hfill\rm(\roman{cnt})\hfill}}{%
\settowidth{\labelwidth}{{\rm(iv)}}\leftmargin=\labelwidth%
\advance\leftmargin by \labelsep\rightmargin=0pt\usecounter{cnt}}}{\end{list}} \makeatletter
\def\mydggeometry{\makeatletter\dg@YGRID=1\dg@XGRID=20\unitlength=0.003pt\makeatother}
\makeatother \theoremstyle{remark}

\numberwithin{equation}{section}

  \DeclareMathOperator{\ad}{ad}
\let\bwdg\bigwedge
\def\bigwedge{{\textstyle\bwdg}}

\begin{document}

\newcommand{\thmref}[1]{Theorem~\ref{#1}}
\newcommand{\secref}[1]{Section~\ref{#1}}
\newcommand{\lemref}[1]{Lemma~\ref{#1}}
\newcommand{\propref}[1]{Proposition~\ref{#1}}
\newcommand{\corref}[1]{Corollary~\ref{#1}}
\newcommand{\remref}[1]{Remark~\ref{#1}}
\newcommand{\defref}[1]{Definition~\ref{#1}}
\newcommand{\er}[1]{(\ref{#1})}
\newcommand{\id}{\operatorname{id}}
\newcommand{\ord}{\operatorname{\emph{ord}}}
\newcommand{\sgn}{\operatorname{sgn}}
\newcommand{\wt}{\operatorname{wt}}
\newcommand{\tensor}{\otimes}
\newcommand{\from}{\leftarrow}
\newcommand{\nc}{\newcommand}
\newcommand{\rnc}{\renewcommand}
\newcommand{\dist}{\operatorname{dist}}
\newcommand{\qbinom}[2]{\genfrac[]{0pt}0{#1}{#2}}
\nc{\cal}{\mathcal} \nc{\goth}{\mathfrak} \rnc{\bold}{\mathbf}
\renewcommand{\frak}{\mathfrak}
\newcommand{\supp}{\operatorname{supp}}
\newcommand{\Irr}{\operatorname{Irr}}
\renewcommand{\Bbb}{\mathbb}
\nc\bomega{{\mbox{\boldmath $\omega$}}} \nc\bpsi{{\mbox{\boldmath $\Psi$}}}
 \nc\balpha{{\mbox{\boldmath $\alpha$}}}
 \nc\bpi{{\mbox{\boldmath $\pi$}}}
\nc\bsigma{{\mbox{\boldmath $\sigma$}}} \nc\bcN{{\mbox{\boldmath $\cal{N}$}}} \nc\bcm{{\mbox{\boldmath $\cal{M}$}}} \nc\bLambda{{\mbox{\boldmath
$\Lambda$}}}

\newcommand{\lie}[1]{\mathfrak{#1}}
\makeatletter
\def\section{\def\@secnumfont{\mdseries}\@startsection{section}{1}%
  \z@{.7\linespacing\@plus\linespacing}{.5\linespacing}%
  {\normalfont\scshape\centering}}
\def\subsection{\def\@secnumfont{\bfseries}\@startsection{subsection}{2}%
  {\parindent}{.5\linespacing\@plus.7\linespacing}{-.5em}%
  {\normalfont\bfseries}}
\makeatother
\def\subl#1{\subsection{}\label{#1}}
 \nc{\Hom}{\operatorname{Hom}}
  \nc{\mode}{\operatorname{mod}}
\nc{\End}{\operatorname{End}} \nc{\wh}[1]{\widehat{#1}} \nc{\Ext}{\operatorname{Ext}} \nc{\textch}{\text{ch}} \nc{\ev}{\operatorname{ev}}
\nc{\Ob}{\operatorname{Ob}} \nc{\soc}{\operatorname{soc}} \nc{\rad}{\operatorname{rad}} \nc{\head}{\operatorname{head}}
\def\Im{\operatorname{Im}}
\def\gr{\operatorname{gr}}
\def\mult{\operatorname{mult}}
\def\Max{\operatorname{Max}}
\def\ann{\operatorname{Ann}}
\def\sym{\operatorname{sym}}
\def\Res{\operatorname{\br^\lambda_A}}
\def\und{\underline}
\def\Lietg{$A_k(\lie{g})(\bsigma,r)$}

 \nc{\Cal}{\cal} \nc{\Xp}[1]{X^+(#1)} \nc{\Xm}[1]{X^-(#1)}
\nc{\on}{\operatorname} \nc{\Z}{{\bold Z}} \nc{\J}{{\cal J}} \nc{\C}{{\bold C}} \nc{\Q}{{\bold Q}}
\renewcommand{\P}{{\cal P}}
\nc{\N}{{\Bbb N}} \nc\boa{\bold a} \nc\bob{\bold b} \nc\boc{\bold c} \nc\bod{\bold d} \nc\boe{\bold e} \nc\bof{\bold f} \nc\bog{\bold g}
\nc\boh{\bold h} \nc\boi{\bold i} \nc\boj{\bold j} \nc\bok{\bold k} \nc\bol{\bold l} \nc\bom{\bold m} \nc\bon{\bold n} \nc\boo{\bold o}
\nc\bop{\bold p} \nc\boq{\bold q} \nc\bor{\bold r} \nc\bos{\bold s} \nc\boT{\bold t} \nc\boF{\bold F} \nc\bou{\bold u} \nc\bov{\bold v}
\nc\bow{\bold w} \nc\boz{\bold z} \nc\boy{\bold y} \nc\ba{\bold A} \nc\bb{\bold B} \nc\bc{\Bbb C} \nc\bd{\bold D} \nc\be{\bold E} \nc\bg{\bold
G} \nc\bh{\bold H} \nc\bi{\bold I} \nc\bj{\bold J} \nc\bk{\bold K} \nc\bl{\bold L} \nc\bm{\bold M} \nc\bn{\Bbb N} \nc\bo{\bold O} \nc\bp{\bold
P} \nc\bq{\Bbb Q} \nc\br{\Bbb R} \nc\bs{\bold S} \nc\bt{\bold T} \nc\bu{\bold U} \nc\bv{\bold V} \nc\bw{\bold W} \nc\bz{\Bbb Z} \nc\bx{\bold
x} \nc\KR{\bold{KR}} \nc\rk{\bold{rk}} \nc\het{\text{ht }}

\nc\toa{\tilde a} \nc\tob{\tilde b} \nc\toc{\tilde c} \nc\tod{\tilde d} \nc\toe{\tilde e} \nc\tof{\tilde f} \nc\tog{\tilde g} \nc\toh{\tilde h}
\nc\toi{\tilde i} \nc\toj{\tilde j} \nc\tok{\tilde k} \nc\tol{\tilde l} \nc\tom{\tilde m} \nc\ton{\tilde n} \nc\too{\tilde o} \nc\toq{\tilde q}
\nc\tor{\tilde r} \nc\tos{\tilde s} \nc\toT{\tilde t} \nc\tou{\tilde u} \nc\tov{\tilde v} \nc\tow{\tilde w} \nc\toz{\tilde z}
\title[Representations of Affine and Toroidal Lie Algebras]{Representations of Affine and Toroidal Lie Algebras\author{Vyjayanthi Chari}}

\begin{abstract} We discuss the category $\cal I$  of level zero integrable representations of loop algebras and their generalizations. The category is not semisimple and so  one is interested in its homological properties. We begin by looking at some approaches which are   used in the study of other well--known non--semisimple categories in the representation theory of Lie algebras. This is done  with a view to seeing if and how far these approaches can be made to work for $\cal I$. In the later sections we focus first  on understanding the irreducible level zero modules and later on certain universal modules, the local and global Weyl modules which in many ways play a role similar to the Verma modules in the BGG--category $\cal O$. In the last section, we discuss the connections with the representation theory of finite--dimensional associative algebras and on some recent work with J. Greenstein.
\end{abstract}

\maketitle

\section*{Introduction} In these notes we discuss aspects of the representation theory of affine  Kac--Moody Lie algebras.  There is a vast literature on this subject and our focus, for the most part, will be on the level zero integrable representations (which include the finite--dimensional  representations) of the loop algebra. At the moment, it appears that this direction is the one most likely to generalize to the case of multi--loop algebras and perhaps to extended affine Lie algebras.

A finite--dimensional  representation of the loop algebra need not be completely reducible and thus one encounters problems of the type  which usually appear in the representation theory of algebraic groups in positive  characteristic or in the BGG--category $\cal O$ for semisimple Lie algebras. However, many of the methods used in the study of the latter subjects are not available for affine Lie algebras. We have tried to  illustrate some of the similarities and the difficulties  throughout the notes.
  Thus, we begin with a brief discussion of  some of the basic ideas, tools and important results in the representation theory of semisimple  Lie algebras and we restrict our attention to ideas that we will use in the later lectures. We then discuss the integrable representations and the BGG--category $\hat{\cal O}$ for affine algebras, the integer form of the universal enveloping algebra of the affine Lie algebra and the positive level integrable modules.  The remaining sections are devoted to integrable level zero modules, the  finite--dimensional representations of the loop algebra and  generalizations of these results to multiloop algebras. We conclude the notes with a section on the connections between the category of graded representations of a maximal parabolic subalgebra of the affine Lie algebra and representations of associative finite--dimensional algebras, highest weight categories and quivers.

We have tried to maintain some level of informality,  reflecting the fact that these are notes of a summer school. This allows for some repetition and topics sometimes appear before they should since they serve as motivation for future lectures.
One is also  guided by one's own interests and necessarily, these notes do not include many important directions in the study of representations of affine Lie algebras: as, for instance, the connections with mathematical physics, number theory, vertex algebras and the relations with the monster group. The interested reader will find, however, that the list of references includes a selection  of books and papers dealing with these aspects of the subject.

\medskip
\noindent
{\bf Acknowledgments.} The organizers of the summer school did a superb job of identifying and attracting excellent graduate students and postdoctoral scholars, and I thank them for giving me the opportunity to be a part of the school. Partial support from the NSF grant DMS-0901253 is also acknowledged.

\section{Simple Lie algebras}\label{C:section 2}
Much of the material in this section can be found in introductory graduate text books \cite{Bo},  \cite{Ca}, \cite{Hu1} on Lie algebras and representation theory. Further references to the literature  may be found throughout the section although these are far from comprehensive.
\subsection{}\label{C:section 2.1} Let $\lie g$ be a finite--dimensional complex simple Lie algebra of rank $n$ with  a fixed  Cartan subalgebra $\lie h$  and let $\Phi$  be the corresponding root system and let $\Pi=\{\alpha_i: i\in I\} $ (where $I=\{1,\cdots,n\}$) be  a  set of simple roots for $\Phi$.  The root lattice $Q$ is the $\bz$--span of the simple roots while $Q^+$  is the $\bn$--span of the simple roots, and the set $\Phi^+=\Phi\cap  Q^+$ is the set of positive roots in $\Phi$. The restriction of  the Killing  form $\kappa:\lie g\times\lie g\to\Bbb C$   to $\lie h\times\lie h$ induces a nondegenerate bilinear form $(\ ,\ )$ on $\lie h^*$ and we let $\{\omega_1,\cdots,\omega_n\}\subset\lie h^*$ be the fundamental weights\index{weight} given by $2(\omega_j,\alpha_i)=\delta_{i,j}(\alpha_i,\alpha_i)$.
 Let $P$ (resp. $P^+$) be the $\Bbb Z$ (resp. $\Bbb N$) span of the $\{\omega_i:1\le i\le n\}$ and note that $Q\subseteq P$. Given $\lambda,\mu\in P$ we say that $\mu\le \lambda$ iff $\lambda-\mu\in Q^+$. Clearly $\le $ is a partial order on $P$. The set $\Phi^+$ has a unique maximal element with respect to this order which is denoted by $\theta$ and is called the highest root of $\Phi^+$. Set $$\rho= \sum_{i=1}^n\omega_i=\frac12\sum_{\alpha\in\Phi^+}\alpha.$$
Let $W$ be the Weyl\index{Weyl group} group of $\lie g$ generated by simple refections $\{s_i: i\in I\}$ and  let~$\ell: W\to \bn$ be the length function which assigns to an element $w\in W$, the length of a  reduced expression for $w$ as a product of simple reflections.

\subsection{}\label{C:section 2.2}  Given $\alpha\in \Phi$, let $\lie g_\alpha$ be the corresponding root space and set $\lie n^\pm=\oplus_{\alpha\in \Phi^+}\lie g_{\pm\alpha}$. Then $\lie n^\pm$ are subalgebras of $\lie g$ and we have an isomorphism of vector spaces \begin{equation}\label{C:gtriang}\lie g=\lie n^-\oplus\lie h\oplus\lie n^+.\end{equation}  For $\alpha\in \Phi^+$, fix elements $x^\pm_\alpha\in\lie g_{\pm\alpha}$ and $h_\alpha\in\lie h$ such that they span a Lie subalgebra of $\lie g$ which is isomorphic to $\lie{sl}_2$, i.e., we have $$[h_\alpha,x^\pm_\alpha]=\pm 2 x^\pm_\alpha,\ \ \ [x^+_\alpha, x^-_\alpha]= h_\alpha, $$ and more generally, assume that the set$$\{x^\pm_\alpha:, \alpha\in\Phi^+\}\cup\{h_{\alpha_i}: i\in I\},$$ is a Chevalley basis for $\lie g$.

\subsection{}\label{C:section 2.3} For any Lie algebra $\lie a$, we let $\bu(\lie a)$ be the universal enveloping algebra of $\lie a$.
Representations of   $\lie a$ are the same as the representations of $\bu(\lie a)$ and since  $\bu(\lie a)$ is  a Hopf algebra, one can define the notions of the trivial representation, the tensor product of representations and the dual of  a representation.
If $\dim\lie a<\infty$, then the algebra $\bu(\lie a)$ is {\em  Noetherian } and this fact plays an important role in the representation theory of the  complex semisimple Lie algebras.  The Poincare--Birkhoff--Witt theorem implies that we have an isomorphism of vector spaces \begin{equation}\label{C:PBW} \bu(\lie g)\cong\bu(\lie n^-)\otimes\bu(\lie h)\otimes \bu(\lie n^+).\end{equation} Another  fact that plays a major role in the representation theory of $\lie g$ is that the center $\cal Z$ of $\bu(\lie g)$ is  large and one has a good theory of central characters, as we shall see below. We recall here the definition of the Casimir\index{Casimir!element} element $\Omega$ in $\cal Z$ and we do this in a way that is adapted for use in the affine case. Let $h_1',\cdots ,h_n'$ be a basis of $\lie h$ which is dual to the basis $h_{\alpha_1},\cdots, h_{\alpha_n}$ with respect to the Killing form of $\lie g$ and set \begin{equation}\label{C:defcasimir}\Omega\index{Casimir!element}=\sum_{i=1}^nh_ih_i'+h_\rho+2\sum_{\alpha\in \Phi^+}x^-_\alpha x^+_\alpha,\end{equation} where $h_\rho\in\lie h$ is defined by  requiring $\rho(h)=\kappa(h,h_\rho)$.

\subsection{}\label{C:section 2.4} Recall that a $\lie g$--module $V$ is said to be a weight module\index{module!weight} if $$V=\bigoplus_{\mu\in\lie h^*} V_\mu,$$ where $V_\mu=\{v\in V: hv=\mu(h)v\ \ h\in\lie h\}$. An element $\mu\in\lie h^*$ is a weight\index{weight} of $V$ if $V_\mu\ne 0$  and we set $\wt(V)=\{\mu\in\lie h^*: V_\mu\ne 0\}$.

The BGG--category $\cal O$\index{Category $\mathcal O$} is the full subcategory of   $\lie g$--modules $M$ satisfying the condition that  $M$  is a finitely generated weight module which is  $\lie n^+$--locally finite\index{locally!finite}  (for all $m\in M$ we have  $\dim\bu(\lie n^+)m<\infty$) as a $\lie n^+$--module. The morphisms in $\cal O$ are just the $\lie g$--module maps. The following is straightforward.
\begin{lem}\label{C:Oprop} The category $\cal O$ is abelian and all objects of $\cal O$ are Noetherian modules for $\bu(\lie g)$. Moreover if $M\in\cal O$ then $$\dim M_\mu<\infty\ \ {\rm{for\ all}}\ \  \mu\in\lie h^*,$$ and there exist finitely many elements  $\mu_1,\cdots,\mu_k\in\lie h^*$ (depending on $M$) such that $$\wt M\subset\bigcup_{s=1}^k\left(\mu_s-Q^+\right). $$\hfill\qedsymbol
\end{lem} The formal  character $\textch(V)$\index{character!formal}  of a weight module $V$, encodes the information on the dimensions of the weight spaces\index{weight!space}, i.e., it is the function $\lie h^*\to\Bbb N$ which sends $\lambda$ to $\dim V_\lambda$. This is usually written as follows.
Given $\lambda\in\lie h^*$ let $e(\lambda): \lie h^*\to \bz$ be defined by $e(\lambda)\lambda=1$, $e(\lambda)\mu=0$ if $\mu\ne \lambda$, then \begin{equation}\label{C:fchar}\textch(V)=\sum_{\mu\in\lie h^*}\dim V_\mu \ e(\mu).\end{equation} Formal characters\index{character!formal} behave well with respect to direct sums and tensor products.

\subsection{}\label{C:section 2.5} An important family of objects in $\cal O$ are the highest weight modules\index{module!highest weight}. A module $M$ is said to be of highest weight  $\lambda$ with highest weight vector $m$ if  $M=\bu(\lie g)m$ and \begin{equation}\label{C:defverma}hm=\lambda(h)m,\quad\ \  x_\alpha^+m=0,\ \ \ h\in\lie h, \ \ \alpha\in\Phi^+.\end{equation}
For $\lambda\in\lie h^*$  let  $M(\lambda)=\bu(\lie g)m_\lambda$ be  the universal highest weight module (also called the  Verma module)\index{module!Verma} with highest weight $\lambda$ and  highest weight vector $m_\lambda$, in other words the equation  \eqref{C:defverma} gives the defining relations of  $m_\lambda$. Using \eqref{C:PBW},
 it is straightforward to show  that  $$\wt M(\lambda)=\lambda-Q^+,\quad \dim M(\lambda)_\lambda=1,\ \ \dim M(\lambda)_\mu<\infty,\quad  \mu\in\lie h^*.$$ In particular, it follows that $M(\lambda)$ has a unique irreducible quotient $V(\lambda)$ and using Lemma \ref{C:Oprop}, we get
 \begin{lem} Any irreducible object in $\cal O$ is isomorphic to $V(\lambda)$ for some $\lambda\in\lie h^*$.\hfill\qedsymbol
 \end{lem}

\subsection{}\label{C:section 2.6} The most obvious family of $\lie g$--modules which are  in $\cal O$ are of course the finite--dimensional ones and  the following  result completely describes the corresponding full subcategory of $\cal O$.

\begin{thm}\hfill

\begin{enumerit}
\item [(i)] Any finite--dimensional $\lie g$--module is isomorphic to a direct sum of irreducible modules, i.e., the full subcategory of $\cal O$ consisting of finite--dimensional modules is semisimple.
    \item[(ii)] For $\lambda\in\lie h^*$, the irreducible module  $V(\lambda)$ is finite--dimensional iff $\lambda\in P^+$ and  the character\index{character} of $V(\lambda)$ is given by the Weyl--character formula, \begin{equation}\label{C_weylc} \left(\sum_{w\in W}(-1)^{\ell(w)} e(w\rho)\right)\textch(V(\lambda))=\sum_{w\in W}(-1)^{\ell(w)}e(w(\lambda+\rho)).
\end{equation}
\item[(iii)]  For $\lambda\in P^+$, the module $V(\lambda)$ is generated by an element $v_\lambda$  with defining relations:$$hv_\lambda=\lambda(h)v_\lambda, \ \ x_{\alpha_i}^+v_\lambda=0,\ \  (x_{\alpha_i}^-)^{\lambda(h_{\alpha_i})+1} v_\lambda=0,$$\ where   $h\in\lie h$ and $1\le i\le n$.

 \end{enumerit}\hfill\qedsymbol
\end{thm}
 We shall say that a $\lie g$--module $V$ is locally finite--dimensional\index{locally!finite-dimensional} if it isomorphic to a sum of finite--dimensional modules. The preceding theorem implies that $V$ is  isomorphic to a direct sum of modules $V(\lambda)$, $\lambda\in P^+$. If in addition, we have $$\dim\Hom_{\lie g}(V(\lambda), V)<\infty\ \ \lambda\in P^+,$$ then we can define a function $P^+\to\bn$ which maps $\lambda$ to $\dim\Hom_{\lie g}(V(\lambda), V)$. We denote this function by $\textch_{\lie g}(V)$ and observe that it  encodes the multiplicity of an irreducible representation in $V$ (whereas $\textch(V)$ encodes the dimension of a weight space). Clearly, we can write,
 \begin{equation}\label{C:gchar} \textch_{\lie g}(V)=\sum_{\lambda\in P^+}\dim\Hom_{\lie g}(V(\lambda), V) e(\lambda),\end{equation} where $e(\lambda)$ is the function (or rather its restriction to $P^+$) defined in Section \ref{C:section 2.4}.

\subsection{}\label{C:section 2.7}  The category $\cal O$ itself is not semisimple, for instance  the Verma--modules  are indecomposable\index{indecomposable} but not always irreducible and we give an example of this, in the case of $\lie{sl}_2$. Let $x,y,h$ be the standard basis of $\lie{sl}_2$. Since $\lie h$ is one--dimensional, we identify the set $P^+$ with  $\bn$ and consider the Verma--module $M(0)$ and its highest weight vector $m_0$. It is straightforward to see that the element $ym_0\in M(0)$  generates a highest weight submodule isomorphic to $M(-2)$ and  one has the following non--split short exact sequence of $\lie{sl_2}$--modules: \begin{equation}\label{C: extverma}  0\to M(-2)\to M(0)\to \bc\to 0.\end{equation}
This lack of semi--simplicity has resulted in the development of many approaches designed to understand the structure of $\cal O$ and we briefly sketch a few directions. The interested reader could refer to \cite{Hu2} and to the references in that book.  In the later sections of these notes, we shall see that some of these approaches  can   or may be adapted to work for   affine Lie algebras while   others will clearly have no parallel for infinite--dimensional Lie algebras.

\subsection{}\label{C:section 2.8} We start with the theory of central characters\index{character!central} as an example of something that is very powerful in the study of $\cal O$ but is essentially not available for affine Lie algebras.  A central character is a homomorphism of algebras $\chi:\cal Z\to \bc$ (where $\cal Z$ is the center of $\bu(\lie g)$) and a $\lie g$--module $M$ is said to admit a central character $\chi$ if:$$ zm=\chi(z)m,\  \ z\in\cal Z,\ \  m\in M.$$ If $z\in\cal Z$, one can prove without too much difficulty that there exists a unique  element $\beta(z)\in\bu(\lie h)$ such that $z-\beta(z)\in\lie n^-\bu(\lie g)\lie n^+$. The corresponding map   $\beta: \cal Z\to\bu(\lie h)$ is an algebra homomorphism and is  called the Harish--Chandra homomorphism\index{Harish-Chandra homomorphism}. Suppose now that $\lambda\in\lie h^*$ and let  $m\in M(\lambda)$, say $m=gm_\lambda$ for some $g\in\bu(\lie g)$. One has $$zm=gzm_\lambda=\lambda(\beta(z))gm_\lambda,$$ where $\lambda\in\lie h^*$ is extended to an algebra homomorphism also denoted $\lambda: \bu(\lie h)\to\bc$. It is  clear that the composite map $\chi_\lambda: \cal Z\to \bc$ sending  $z\to\lambda(\beta(z))$ is an algebra homomorphism and hence $M(\lambda)$ admits a central character\index{character!central}.

It is now  natural to ask if the $\chi_\lambda$ are all distinct and if any homomorphism from $\cal Z\to\bc$ is of the form $\chi_\lambda$ for some $\lambda\in\lie h^*$. The answer to the first question is quite clearly no. Consider the example given in \eqref{C: extverma}. Since $M(0)$ admits a central character\index{character!central} any submodule or quotient of $M(0)$ also has the same  central character and hence we get in this case that $\chi_0=\chi_{-2}$. It is a theorem of Harish--Chandra, that for $\lambda,\mu\in\lie h^*$ : \begin{equation}\label{C:dotweyl}\chi_\lambda=\chi_\mu \ \iff \ \ \lambda+\rho =w(\mu+\rho)\ \ {\rm{for\ some}}\ \ w\in W,\end{equation} and also that the answer to the second question above is yes. As a   consequence of this theorem one can show that:
\begin{thm} Any module $M\in\Ob\cal O$ is Artinian\index{module!Artinian}. In particular, $M$  has a Jordan--H\"older series and can be written uniquely (up to isomorphism and re--indexing) as a direct sum of indecomposable\index{indecomposable} modules.\hfill\qedsymbol
\end{thm}
For $M\in\cal O$ and $\mu\in\lie h^*$ denote by $[M:V(\mu)]$ the multiplicity of $V(\mu)$ in the Jordan--H\"older series of $M$. The problem of determining the multiplicities of the composition factors\index{composition factor} of the Verma module is  a very difficult one and the answer   is given by  the  famous Kazhdan--Lusztig conjecture  \cite{KL} in terms of  the Kazhdan--Lusztig  (KL) polynomials. The conjecture is now a theorem due to Beilinson--Bernstein \cite{BB}  and independently Brylinski--Kashiwara \cite{BK}.  The KL--polynomials are defined recursively  and explicit computations using the definition are formidable.  There is extensive literature devoted  to understanding the combinatorics of the Kazhdan--Lusztig polynomials and to developing programs to compute the polynomials.

\subsection{} \label{C:section 2.9} It is natural to consider the homological properties of $\cal O$  by asking questions such as: does $\cal O$  contain enough projectives? In  fact it does and  one can determine all indecomposable\index{indecomposable} projective objects   in $\cal O$.  These are again indexed by elements of $\lie h^*$ and we denote by $P(\lambda)$ the corresponding indecomposable\index{indecomposable} projective module\index{module!projective}. The module $P(\lambda)$ has  a  Verma  flag which is  a filtration in which the subsequent quotients are Verma--modules. As in the case of Jordan--H\"older series of a module, the length of any two  filtrations of $P(\lambda)$ by Verma modules is the same and so is  the multiplicity of a Verma module $M(\mu)$ in any two filtrations.

 And this allows us now to state the fundamental result proved by Bernstein--Gelfand and Gelfand:
  \begin{thm}\label{C:bgg} For  $\lambda, \mu\in\lie h^*$ we have $$[P(\lambda): M(\mu)]=[M(\mu): V(\lambda)],$$ where $[P(\lambda): M(\mu)]$ is the multiplicity of $M(\mu)$ in a Verma flag\index{Verma flag} of $P(\lambda)$ and $[M(\mu): V(\lambda)]$ is the multiplicity of $V(\mu)$ in  the Jordan--H\"older series of $M(\mu)$.\hfill\qedsymbol\end{thm} Much of the work on $\cal O$ including the Kazhdan--Lusztig conjecture  was stimulated by this result.

 \subsection{}\label{C:section 2.10} One could try to understand $\cal O$ by looking  at smaller, more manageable subcategories of $\cal O$ and this is done as follows. It is not true that the center always acts on  objects (even indecomposable\index{indecomposable} ones)  of $\cal O$ via a central character\index{character!central}. However, since $\cal Z$ is commutative and we are working over the complex numbers it is true that any $M\in\Ob\cal O$ can be written as a direct sum of generalized eigenspaces for the action of $\cal Z$. Hence if $M$ is indecomposable\index{indecomposable} then there exists an
algebra homomorphism $\chi: \cal Z\to \bc$ such that $$M=M^\chi=\{m\in M: (z-\chi(z))^k m=0\ \ {\rm for \ some}\ k =k(m)\in\bn\}.$$ Moreover in this case, if $V(\mu)$ occurs in a Jordan--H\"older series for $M$ then $\chi=\chi_\mu$.  This motivates the following definition. Let $\cal O_\chi$ be the full subcategory of $M$ consisting of $M\in\Ob\cal O$ such that $M=M^\chi$. Using \eqref{C:dotweyl} we see that $\cal O_\chi$ has only finitely many simple objects $V(\mu)$ with $\chi_\mu=\chi$  and it suffices to understand  $\cal O_\chi$.  The endomorphism algebra of the projective generator of $\cal O_\chi$ is a finite--dimensional associative algebra $\cal A_\chi$ and its (left)--module category is equivalent to $\cal O_\chi$. The category $\cal O_\chi$ fits into the axiomatic framework of highest weight categories developed by Cline, Parshall and Scott, \cite{CPS} and  $\cal A_\chi$ is quasi--hereditary. The algebra $\cal A_\chi$ is given by a quiver with relations although it is in general very hard to compute these relations, see however \cite{St} where an algorithm to compute these relations is given. Also, under suitable conditions  \cite{BGS}, \cite{So}, the algebra $\cal A_\chi$ has a Koszul grading and hence the homological properties of $\cal O_\chi$ can also be understood from this viewpoint.

\subsection{}\label{C:charpg}\label{C:section 2.11} Before we leave the realm of finite--dimensional simple Lie algebras,  we discuss the Kostant $\Bbb Z$--form\index{Kostant $\Bbb Z$--form} of $\bu(\lie g)$ and the passage to positive characteristic. Given any associative algebra $A$ over $\Bbb C$ and elements $a\in A$ and $r\in\Bbb N$, the $r^{th}$--divided power of $a$ is the element $a^{(r)}=a^r/r!$.
 Let $\bu_{\Bbb Z}(\lie g)$ be the $\Bbb Z$--subalgebra of $\bu(\lie g)$ generated by the elements $(x_\alpha^\pm)^{(r)}$, $\alpha\in \Phi^+$, $r\in\Bbb N$ and define $\bu_{\Bbb Z}(\lie n^\pm)$ in the obvious way. To define the analog of the Cartan subalgebra, one does not do the obvious and take the divided powers of the $h_i$, $1\le i\le n$. Instead, one is guided by the following formula which rewrites the product $(x^+_\alpha)^{(s)}(x^-_\alpha)^{(r)}$ in Poincare--Birkhoff--Witt order as in \eqref{C:PBW},
 \begin{equation}\label{C:zform}(x^+_\alpha)^{(s)}(x^-_\alpha)^{(r)}=\sum_{m=0}^{\min(r,s)}(x^-_\alpha)^{(r-m)}\binom{h_\alpha-r-s+2m}{m}(x^+_\alpha)^{(s-m)},
 \end{equation}
 where for any $m\in\Bbb N$ and $h\in\lie h$, we set $$\binom{h}{m}=\frac{h(h-1)\cdots (h-m+1)}{m!}.$$ We then define $\bu_{\Bbb Z}(\lie h)$ to be  the $\Bbb Z$--subalgebra generated by the elements $\binom{h_{\alpha_i}}{m}$, $1\le i\le n$, $m\in\Bbb N$. It is a theorem of Kostant that  $\bu_{\Bbb Z}(\lie g)$ has a PBW--basis, by which we mean  a $\bz$--basis of $\bu_\bz(\lie h)$ consisting of all ordered monomials from the set $$\{(x^\pm_\alpha)^{(r)}, \binom{h_i}{r}: \alpha\in \Phi^+,\ \ i\in I, r\in\Bbb N\}.$$ This means that $\bu_{\Bbb Z}(\lie g)$ is a $\Bbb Z$--lattice\index{lattice} in $\bu(\lie g)$ and moreover one can prove that if $\lambda\in P^+$, and we set $$V_{\Bbb Z}(\lambda)=\bu_{\Bbb Z}(\lie g)v_\lambda,$$ then $V_{\Bbb Z}(\lambda)$ is also a $\Bbb Z$--lattice in $V(\lambda)$. Let $\Bbb F$ be a field of characteristic $p$ and regard $\Bbb F$ as a left module for $\Bbb Z$.  Set $$\bu_{\Bbb F}(\lie g)=\bu_{\Bbb Z}(\lie g)\otimes_{\Bbb Z}\Bbb F,\ \ \ \  V_{\Bbb F}(\lambda)=V_{\Bbb Z}(\lambda)\otimes_{\Bbb Z}\Bbb F.$$ Clearly, $\bu_{\Bbb F}(\lie g)$ is an associative algebra over $\Bbb F$ (called the hyperalgebra\index{hyperalgebra} of $\lie g$ over $\Bbb F$) and $V_{\Bbb F}(\lambda)$ is an indecomposable\index{indecomposable}  (usually reducible)  module for this algebra, and is called the Weyl module\index{module!Weyl} of weight $\lambda$. There is extensive literature on the  study of finite--dimensional representations of the hyperalgebra\index{hyperalgebra}.  One can talk about weight\index{module!weight} modules  (whose definition requires a natural  modification) and characters\index{character} of these modules and we summarize the relevant points for our  next lectures in the following theorem.
  \begin{thm} \label{C:domfin}\begin{enumerit}
  \item [(i)]  Given $\lambda\in P^+$, the module $V_{\Bbb F}(\lambda)$ has a unique irreducible quotient and any irreducible finite--dimensional module of $\bu_{\Bbb F}(\lie g)$ is isomorphic to such a quotient for some $\lambda\in P^+$.
      \item[(ii)] The character\index{character} of  $V_{\Bbb F}(\lambda)$ is given by the Weyl character\index{character!formula} formula.
      \item[(iii)]  The module $V_{\Bbb F}(\lambda)$ is generated by the element $v_\lambda\otimes 1$ and the defining relations are:
      $$(x^+_\alpha)^{(r)}(v_\lambda\otimes 1)=0,\ \  \binom{h_i}{r}(v_\lambda\otimes 1)=\binom{\lambda(h_i)}{r}(v_\lambda\otimes 1),\ \ (x^-_\alpha)^{(s)}(v_\lambda\otimes 1)=0,$$ for all $\alpha\in \Phi^+$, $i\in I$,  $r,s\in\Bbb N$ with $s\geq\lambda(h_i)+1$.
      \end{enumerit}\hfill\qedsymbol\end{thm}
      Note that the fact that the modules $V_{\Bbb F}(\lambda)$ are indecomposable\index{indecomposable} but not irreducible implies  that the category of finite--dimensional representations of the hyperalgebra\index{hyperalgebra} is also not semisimple and many of the problems and methods discussed for the category $\cal O$ are also studied and used in this case.

      \section{Affine Lie algebras}\label{C:section 3}  We now turn our attention to the case of the best understood Kac--Moody Lie algebras: the affine Lie algebras\index{affine Lie algebra}. There are two definitions available for these algebras: one is via a finite set of generators and relations and the other is very explicit as we shall see below.   It  is the interplay between these two definitions that makes the study of these algebras and their representations tractable. The explicit realization is what we need for  our purposes and so we refer the reader to \cite{Ca}, \cite{Ka2} for the presentation via Chevalley generators and relations. Most of the results  discussed in this section can also be found in those two books. We use freely the notation of the first section.

      \subsection{}\label{C:section 3.1}   Let $t$ be an indeterminate and $\bc[t,t^{-1}]$  the algebra of Laurent polynomials in $t$ with the coefficients in $\bc$. For any complex Lie algebra $\lie a$, the  loop algebra\index{loop algebra} $L(\lie a)$  is defined by, $$L(\lie a)=\lie a\otimes\bc[t, t^{-1}]$$ with Lie bracket, $$ [a\otimes f(t), b\otimes g(t)]=[a,b]\otimes f(t)g(t),\ \ a,b\in \lie a,\ \ f,g\in\bc[t,t^{-1}].$$ The algebra $\lie a$ is identified with the subalgebra $\lie a\otimes 1$ and we shall assume this from now on without further mention.
Given  an automorphism $\sigma$ of $\lie a$ of order $m$ and $\zeta$  a fixed primitive $m^{th}$ root of unity, we can write,   $$\lie a=\bigoplus_{s=0}^{m-1}\lie a_s,\ \qquad \ \lie a_s=\{a\in\lie a: \sigma(a)=\zeta^s a\},$$ and we set $$L(\lie a,\sigma, m)= \bigoplus_{s=0}^{m-1}\lie a_s\otimes t^s\bc[t^m,t^{-m}].$$ It is a simple matter to check that $L(\lie a, \sigma, m)$ is a Lie subalgebra of $L(\lie a)$ and it is a proper subalgebra if and only if  $\sigma$ is not the identity automorphism of $\lie a$, and in this case, we call $L(\lie a, \sigma, m)$ a twisted loop algebra\index{loop algebra!twisted}.

     \subsection{}\label{C:section 3.2}   The  Lie algebra   $L(\lie g)$ admits a unique (up to an  isomorphism)  non--trivial central extension which we denote by $\tilde L(\lie g)$. It can be constructed as   follows: set $$\tilde L(\lie g)= L(\lie g)\oplus\bc c$$ and define the Lie bracket by requiring $c$ to be central and by setting $$ [x\otimes t^p , y\otimes t^r] =[x,y]\otimes t^{p
     +r}+ p\delta_{p+r,0}\kappa(x,y)c,\ \ x,y\in\lie g,\ \ p,r\in\bz,$$ where $\delta_{p+r,0}$ is the Kronecker $\delta$ symbol and recall that $\kappa$ is the Killing form of $\lie g$. The derivation $d=td/dt$ of $\bc[t,t^{-1}]$ acts on $L(\lie g)$ and on $\tilde L({\lie g})$ by  $$d(x\otimes t^r)= rx\otimes t^r,\ \ \ [d,c]=0.$$ The semi--direct product Lie algebra $$\hat{L}(\lie g)= \tilde L(\lie g)\oplus \bc d,$$ is called the  untwisted affine Lie algebra and  $\lie g$ is  naturally isomorphic to a subalgebra of $\hat{L}(\lie g)$. If $\sigma$ is a Dynkin diagram automorphism\index{diagram automorphism} of $\lie g$, then  $\tilde L(\lie g, \sigma, m)$ and $\hat L(\lie g, \sigma, m)$ are defined in the obvious way and are Lie subalgebras of $\tilde L(\lie g)$ and $\hat{L}(\lie g)$ respectively. The algebra $\hat L(\lie g, \sigma, m)$ is called a twisted affine Lie algebra\index{affine Lie algebra!twisted}.

     {\em For the most part, we will deal only with the untwisted  affine Lie algebras and we will refer to these simply as affine Lie algebras. We  will usually  indicate if the  results we discuss are known for  for the twisted affine algebras.}

      \subsection{}\label{C:section 3.3}  The affine Lie algebra  comes equipped with a naturally defined root system and a set of simple roots. Set $\hat{\lie h}=\lie h\oplus\bc c\oplus\bc d$ and define $\delta\in \hat{\lie h} ^*$  by  $\delta(d)=1$, $\delta(\lie h\oplus\bc c)=0$.  Then $\hat{\lie h}$ is abelian and given $\mu\in\lie h^*$ we extend $\mu$ to an element of $\hat{\lie h}^*$ by requiring $$\mu( d)=0,\ \qquad  \mu(c)=0,$$ and by abuse of notation, we continue to denote this element by $\mu$.  The  adjoint action of $\hat{\lie h}$ on $\hat L(\lie g)$ is  semisimple and  the non--zero  eigenvalues (the set of roots)  are, $$\hat{\Phi}=\{\alpha+r\delta: \alpha\in \Phi, r\in\bz\}\bigsqcup\{r\delta: r\in\bz, \   \   r\ne 0\}.$$ The corresponding eigenspaces (root spaces) are $\lie g_\alpha\otimes t^s$, $\lie h\otimes t^r$ where $\alpha\in \Phi$, $r,s\in\bz$ and $r\ne 0$. Notice that the eigenspaces corresponding to $r\delta$ are of dimension $\dim\lie h$ and hence bigger than one in general.  Setting $$\hat\Pi=\{\alpha_i: 0\le i\le n\},\ \ \alpha_0= -\theta+\delta,$$  we find that any element of $\hat{\Phi}$ can be written uniquely as an integer linear combination of elements of $\hat{\Pi}$ where the coefficients are all either non--negative or non--positive. Elements of the set $\hat\Pi$ are called simple roots and we have $$\hat\Phi^+=\{\pm\alpha+r\delta: r\in\bn_+,\ \ \alpha\in \Phi\}\quad \sqcup \ \Phi^+\quad \sqcup\ \{r\delta: \ r\in\bn_+\}.$$ The subalgebras $\hat{\lie n}^\pm$ are defined in the natural way and one has  triangular decompositions similar to the ones given in \eqref{C:gtriang} and \eqref{C:PBW} for $\lie g$. Set $$x_\theta^\mp\otimes t^{\pm 1} =x^\pm_{\alpha_0},\qquad  h_0= c-h_\theta.$$  The  elements $x^\pm_{\alpha_i}$, $h_{\alpha_i}$ $0\le i\le n$  are called the Chevalley generators of  $\hat L(\lie g)$ and an abstract presentation of the Lie algebra in terms of these generators can be given. The fundamental weights\index{weight} $\Lambda_i$, $0\le i\le n$ are elements of $\hat{\lie h}^*$ defined by $\Lambda_i(h_{\alpha_j})=\delta_{i,j}$.
The notion of root lattice $\hat Q$, weight\index{weight!lattice} lattice $\hat P$, weight modules, weight spaces, highest weight modules and so on, have obvious analogs in the affine case, one just replaces $\lie h$ by $\hat{\lie h}$ and so on.  The Killing form of $\hat L(\lie g)$ cannot be defined in the usual way since $\hat L(\lie g)$ is infinite--dimensional. However, the formulae
\begin{gather*}< x\otimes t^r, y\otimes t^s>=\kappa(x,y)\delta_{r+s, 0},\ \ <x\otimes  t^r, \bc c\oplus\bc  d>=0,\\ <c,c>=<d,d>=0,\ \ <c,d>=1,\end{gather*} define an invariant symmetric nondegenerate bilinear form on $\hat L(\lie g)$. The affine Weyl group\index{Weyl group!affine} $\hat W$  is the subgroup of ${\rm{Aut}}(\lie h^*\oplus\bc\delta)$ generated by the simple reflections $\{s_i: i=0,\cdots ,n\}$ and is an infinite Coxeter group. As for finite--dimensional simple Lie algebras, the length of an element $w$ in $\hat W$ is just the number of positive roots which are turned negative by $w$: in particular (unlike in the situation for finite--dimensional simple Lie algebras) there does not exist an element in $\hat W$ which maps $\hat\Phi^+$ to $\hat\Phi^-$.

        The many parallels with the structure of $\lie g$ suggest that the study of the representation theory of $\hat L(\lie g)$  should proceed in the same way as the study of representations of $\lie g$. However, the differences noted above,
         together with the observation that $\bu(\hat L(\lie g))$ is not Noetherian  means that there are significant  differences  and  difficulties   in the representation theory of $\hat L(\lie g)$ and we shall see some of these  in the rest of the notes.

    \subsection{}\label{C:section 3.4}    Before defining the category $\hat{\cal O}$,  we note the following result proved in \cite{CI}:
      \begin{prop} The center of $\bu(\hat{L}(\lie g))$ is the polynomial algebra generated by the center $\bc c$ of $\hat L(\lie g)$.\hfill\qedsymbol
      \end{prop}
      This proposition means that unlike Section 2,  there is no good theory of central characters\index{character!central} available for the study of representations of  affine Lie algebras.  However, many of the results discussed for $\cal O$ can be proved for $\hat{\cal O}$ and this is essentially because it is  possible to define the analog of  the Casimir operator\index{Casimir!operator}. Recall that the  definition given in  \eqref{C:defcasimir} of the Casimir, only used dual bases with respect to the Killing form of $\lie g$. This formula can be used verbatim for $\hat L(\lie g)$ by using the symmetric form defined above.  However the sum is infinite and the Casimir operator lives in some completion of  $\bu(\hat L(\lie g))$ and this means that it does not operate on all $\hat L(\lie g)$--modules. But it does act on modules in the category $\hat{\cal O}$ whose definition we now recall.

    \subsection{}\label{C:section 3.5}     We say that an $\hat L(\lie g)$--module  $M$ is an object of $\hat{\cal O}$  if it is a weight\index{module!weight} module, i.e., $\hat{\lie h}$ acts semi--simply,  $ M=\oplus_{\mu\in\hat{\lie h}^*} M_\mu$, and \begin{equation}\label{C:cone} \dim M_\mu<\infty,\ \ \wt M\subset \bigcup_{s=1}^p\mu_s-\hat Q^+,\end{equation} for some $p\in\bn$ and elements $\mu_1,\cdots,\mu_p\in\hat{\lie h}^*$. This definition differs from that of $\cal O$, in two ways.   The first is that one does not require the module to be finitely generated.  This is because one wants $\cal O$ to be closed under taking submodules and quotients and since $\bu(\hat L(\lie g))$ is not Noetherian a submodule of a finitely generated module need not be finitely generated. The condition that $\hat{\lie n}^+$ acts locally finitely on $M$ does not have all the implications it does in the category $\cal O$  (see Lemma \ref{C:Oprop}) and  hence is replaced by the stronger conditions in  \eqref{C:cone}. Moreover, notice also that the second condition in \eqref{C:cone} implies that if $\mu\in\wt M$, then $\mu+\beta\in\wt M$ only for finitely many $\beta\in\hat\Phi^+$. It is now immediate that the Casimir\index{Casimir!operator} operator acts in a well--defined way on every $M\in\hat{\cal O}$ and moreover one can prove that it   commutes with the action of $\hat L(\lie g)$.

       \subsection{}\label{C:section 3.6}   For $\lambda\in\hat{\lie h}^*$ we let $\hat M(\lambda)$ be the Verma module with highest weight\index{weight} $\lambda$ defined (with obvious modifications) in the same way as  the Verma modules for simple Lie algebras given in Section \ref{C:section 2.5}. It has a unique irreducible quotient $\hat V(\lambda)$ and any irreducible module in $\hat{\cal O}$ is isomorphic to one of these and these results are proved in the same way as those for simple Lie algebras.

        The following is a fairly simple exercise.
       \begin{lem} Let $M\in\hat{\cal O}$ and assume that $\dim M<\infty$. Then $(L(\lie g)\oplus \bc c)M=0.$\hfill\qedsymbol
       \end{lem}
      \noindent  In other words, there are no interesting finite dimensional modules in $\hat{\cal O}$ and in view of Theorem \ref{C:section 2.6}  this raises the natural question,   are irreducible  modules associated to $\lambda\in\hat P^+$ distinguished in some way? The answer is yes and this leads us to the definition of an integrable module\index{module!integrable} for an   affine Lie  algebra.
     \subsection{}\label{C:section 3.7}   A weight\index{module!weight} module  $\hat V$  for $\hat L(\lie g)$ is said to be  integrable\index{module!integrable} if the elements $x^\pm_{\alpha_i}$, $0\le i\le n$  act locally nilpotently\index{locally!nilpotent} on $\hat V$.
      Integrable\index{module!integrable} modules share many of the properties of finite--dimensional $\lie g$ modules. For instance a standard $\lie{sl}_2$--argument shows that if $V$ is integrable then $\wt V\subset \hat P$ and also that $\wt V$ is   $\hat W$--invariant. (Of course, one could have also defined integrable modules for $\lie g$ but the definition is not interesting since one can prove that any integrable $\lie g$--module is a sum of finite--dimensional $\lie g$--modules). Let $\textch\hat V(\lambda)$ be the formal character of $\hat V(\lambda)$ as in \eqref{C:fchar}. The following theorem explains the analogy between finite--dimensional modules for $\lie g$ and integrable modules  in $\hat{\cal O}$.
   \begin{thm} For  $\lambda\in\hat P^+$, the module $\hat V(\lambda)$ is integrable\index{module!integrable} and the full subcategory of $\hat{\cal O}$ consisting of integrable modules is completely reducible. Finally, \begin{equation}\label{C:weyla} \left(\sum_{w\in \hat W}(-1)^{\ell(w)} e(w\rho)\right)\textch(\hat V(\lambda))=\sum_{w\in\hat  W}(-1)^{\ell(w)}e(w(\lambda+\rho)).\end{equation}\hfill\qedsymbol

   \end{thm}
  The formula \eqref{C:weyla} is called the Weyl--Kac character formula\index{character!formula} and the  Casimir operator\index{Casimir!operator} $\hat{\Omega}$ is an adequate substitute for the role played by the center in proving the Weyl character formula for a simple  Lie algebra. The sums that appear in the Weyl--Kac character formula\index{character!formula} are infinite.
This formula can be specialized in many ways and leads to interesting number--theoretic identities, including combinatorial identities of Macdonald and the classical Jacobi triple product identity.
\subsection{}\label{C:section 3.8}  Suppose that  $V$ is a weight\index{module!weight} module for $\hat L(\lie g)$,  in which case we have  $$V= \oplus_{a\in\Bbb C} V^a,\ \ \ V^a=\{v\in V: cv=av\},$$and $V^a$ is a $\hat L(\lie g)$--submodule. We say that $V^a$ is a module of
level $a$\index{level}.
\begin{lem} For  $\lambda\in\hat{\lie h}^*$, the module $\hat V(\lambda)$ has level $\lambda(h_\theta+h_0)$. Hence every irreducible integrable\index{module!integrable} module in $\hat{\cal O}$ has non--negative integer level. The only level zero irreducible modules in $\hat{\cal O}$ are one--dimensional and correspond to taking $\lambda$ to be a scalar multiple of $\delta$.\hfill\qedsymbol

\end{lem}
Thus all the affine Lie algebras have a canonical integrable\index{module!integrable} representation of level one corresponding to the weight\index{weight} $\Lambda_0$ and this is sometimes called the basic representation\index{basic representation} of $\hat L(\lie g)$. It is the simplest integrable representation in $\hat{\cal O}$ and the character formula given in Section \ref{C:section 3.7} can be made explicit.  The theory of vertex algebras,  the relationship with physics and the connections with the monster group (see \cite{Bor}, \cite{FK}, \cite{FLM}, \cite{Ka2}, for instance) all have their roots in  the effort to understand and explicitly construct the basic representation\index{basic representation} of the affine Lie algebra.

 To conclude this section, we note that many of the problems studied for $\cal O$ have also been studied for $\hat{\cal O}$ and  one does have a block\index{block} decomposition \cite{DGK} and a  Kazhdan--Lusztig theory (see \cite{Hu2} for references to the literature on these topics). But the fact that the category $\hat{\cal O}$ is neither   Noetherian  nor Artinian does make things much more complicated. In the remaining sections our focus will be only on the integrable\index{module!integrable} representations of affine Lie algebras.  Finally, note that the results of this section are also known for the twisted affine algebras.

\section{Affine Lie algebras integrable representations and integral forms}\label{C:section 4} The  reader has noticed by now that there are integrable\index{module!integrable} representations of $\hat L(\lie g)$ which are not in $\hat{\cal O}$.  The most obvious representation of a Lie algebra is the adjoint representation, and it is a simple matter to see from Section \ref{C:section 3.2} that the adjoint representation of $\hat L(\lie g)$ is integrable but not in $\hat{\cal O}$. This raises the problem of classifying the irreducible integrable representations of affine Lie algebras, rather than just the ones in $\hat{\cal O}$. We  address this problem  in the first part of this section.

   The adjoint representation is an example of a level zero representation which is indecomposable\index{indecomposable} and reducible, since the center of $\hat L(\lie g)$ is a proper non--split submodule under the adjoint action.
 This  shows that the category of level zero integrable\index{module!integrable} representations of $\hat L(\lie g)$ is not semisimple and hence should have interesting homological properties.  To study these properties and to pursue  the directions of study which have been fruitful in  the case of simple Lie algebras, it is helpful to understand (following H. Garland) the integral form\index{integral form} of the universal enveloping algebra of $\hat L(\lie g)$. Garland also proved that the  integrable\index{module!integrable} modules in $\hat{\cal O}$ have an integral form, but relatively  little is known about the corresponding representations of the hyper-algebra, except in the case of the basic representation\index{basic representation} which we discuss briefly.

  \subsection{}\label{C:section 4.1} We begin by making some comments about duality.  Let $w_0\in W$ be the longest element of the Weyl group\index{Weyl group} of $\lie g$. It is well--known that the dual of an irreducible  finite--dimensional representation of $\lie g$ with highest weight\index{weight} $\lambda$  is  a highest weight representation with highest weight\index{weight} $-w_0\lambda$. Let us now consider the situation of irreducible integrable\index{module!integrable} modules in $\hat{\cal O}$. The first difficulty one encounters, is that the module $\hat V(\lambda)$ is infinite--dimensional and the dual space is too big. The weight\index{weight!space} spaces, however, are finite--dimensional and so one works instead with the restricted dual\index{restricted dual} $$\hat V(\lambda)^\#=\bigoplus_{\mu\in\hat P^+} V(\lambda)_\mu^*,$$ which is an integrable\index{module!integrable}  $\hat L(\lie g)$--module. Since $\hat W$ is an infinite group, one does not have the analog of the longest element: there does not exist $w\in\hat W$ such that $w\hat{\Pi}=-\hat{\Pi}$ and this means that $\hat V(\lambda)^\#$ is not a highest weight\index{module!highest weight} module. However, it is easily seen that it is a lowest weight module, which is defined in the obvious way, by replacing $\hat{\lie n}^+$ by $\hat{\lie n}^-$. The lowest weight space is $V(\lambda)_\lambda^*$ and the lowest weight is $-\lambda$.  But apart from this modification, it is clear that understanding highest weight modules is the same as understanding lowest weight modules.

  \subsection{}\label{C:section 4.2}\label{C:loops} The category of all integrable\index{module!integrable} modules for $\hat L(\lie g)$ is very big, even if we restrict to a fixed level. To make the study more manageable we let $\cal I$ be the category whose objects are $\hat L(\lie g)$--modules with finite--dimensional weight\index{weight!space} spaces. Our goal is to classify the irreducible objects in $\cal I$, and in the previous section we saw how to construct irreducible objects of $\cal I$ which are also in $\hat{\cal O}$. But so far we have only encountered one interesting example (see Lemma \ref{C:section 3.8}) of a level zero representation, the adjoint, and that was reducible. Our first task then, is to construct natural examples of irreducible representations of level zero.

   Given a representation $V$  of $\lie g$, set $L(V)=V\otimes \bc[t,t^{-1}]$. For $a\in\bc^\times$, one can define the structure of a $\hat L(\lie g)$--module on $L(V)$ by:
  \begin{equation}\label{C:loopv} (x\otimes t^r)(v\otimes t^s)= a^rxv\otimes t^{r+s},\ \ \ c (L(V))=0,\ \ \ d(v\otimes t^r)=rv\otimes t^r.
  \end{equation}
  We denote this module by $L(\ev_a^*V)$ and we will say more about this notation later.
 This construction can be generalized further as follows. Suppose that we are given representations $\{V_s:1\le s\le k\}$ of $\lie g$ and $k$ non--zero complex numbers $a_1,\cdots, a_k$. The space $L(V_1\otimes\cdots\otimes V_k)$ has a $\hat L(\lie g)$--module structure defined by
  \begin{gather*} (x\otimes t^r)(v_1\otimes\cdots\otimes v_k\otimes t^m) =\sum_{s=1}^k(v_1\otimes\cdots\otimes v_{s-1}\otimes xv_s\otimes v_{s+1}\otimes\cdots\otimes v_k\otimes a_s^rt^{r+m}), \\ \\ c L(V_1\otimes\cdots\otimes V_k)=0,\\ \\  d(v_1\otimes\cdots\otimes v_k\otimes t^m) = m (v_1\otimes\cdots\otimes v_k\otimes t^m) , \end{gather*}  and we denote this module by $L(\ev_{a_1}^*V_1\otimes\cdots\otimes \ev_{a_k}^* V_k)$. Modules of these kind are sometimes called loop modules\index{module!loop}.
The following result was proved in \cite{CP}.
  \begin{prop}\label{C:loopirr} Let $k\in\bn$, and  for $1\le s\le k$ let $\lambda_s\in P^+$ and  $a_s\in\bc^\times$  and assume that $a_s\ne a_m$ if $s\ne m$. There exists  $r\ge 1$ such  that $$\sum_{s=1}^k\lambda_s(h) a_s^j=0,\ \  1\le j <r,\; \text{for all }h\in\lie h.$$ For $1\le j\le r$, the submodule $$L^j(\ev_{a_1}^*V(\lambda_1)\otimes\cdots\otimes \ev_{a_k}^* V(\lambda_k))= \bu(\hat L(\lie g))( v_{\lambda_1}\otimes \cdots\otimes v_{\lambda_k}\otimes t^j)$$ is irreducible and  we have an isomorphism of $\hat L(\lie g)$--modules, $$L(\ev_{a_1}^*V(\lambda_1)\otimes\cdots\otimes \ev_{a_k}^* V(\lambda_k))\cong\bigoplus_{s=1}^rL^s(\ev_{a_1}^*V(\lambda_1)\otimes\cdots\otimes \ev_{a_k}^* V(\lambda_k)). $$\hfill\qedsymbol
   \end{prop} The statement in \cite{CP} is  more precise and gives the value of $r$ for a fixed set of  $\lambda_s$, $a_s$, $1\le s\le k$.  As an example, one has    $r=2$ when $\lambda_1=\lambda_2$ and $a_1=-a_2$.

\subsection{}\label{C:section 4.3}   The following  is an amalgamation of results  proved in \cite{Cinv} and  \cite{CP} and achieves our goal of classifying the irreducible objects of $\cal I$. Recall that if $V$ is integrable\index{module!integrable}, then $\wt V\subset\hat P$ and hence the center $c=h_0+h_\theta$ acts on $V$ with integer eigenvalues.
  \begin{thm}\label{C:inv} Let $V\in\Ob\cal I$ be irreducible and let  $k\in\bz$ be such that  $c v=k v$ for all $v\in V$.
 \begin{enumerit}
 \item[(i)] If $k>0$ then $V\in\hat{\cal O}$ and hence $V\cong\hat{V}(\lambda)$ for some $\lambda\in\hat P^+$.
 \item[(ii)] If $k<0$ then $V\cong\hat{V}(\lambda)^\#$ for some $\lambda\in\hat P^+$.
 \item[(iii)] If $k=0$, then  $$V\cong L^s(\ev_{a_1}^*V(\lambda_1)\otimes\cdots\otimes \ev_{a_k}^* V(\lambda_k)),$$ for some     $\lambda_1,\cdots,\lambda_k\in P^+$, $a_1,\cdots, a_k\in\bc^\times$ and $s$  as in Proposition \ref{C:loopirr}. \end{enumerit}
\hfill\qedsymbol \end{thm}
The twisted version of this theorem is known and may be found in \cite{CPtwisted}.

 We now make some remarks  which will allow us to limit our  focus in the next sections on finite--dimensional representations of $L(\lie g)$. The first observation is that  to produce a representation of $\hat L(\lie g)$ from a representation of $\lie g$  one has to tensor with $\bc[t,t^{-1}]$  only to keep track of the grading that $d$ defines on $\hat L(\lie g)$.  But one can get around this problem by defining a functor $\bl$ from the category of $L(\lie g)$--modules to  $\hat L(\lie g)$--modules, which is given on objects by, \begin{gather*} \bl(V) =V\otimes \bc[t,t^{-1}], \\ (x\otimes t^r)(v\otimes t^s)= (x\otimes t^rv)\otimes t^{r+s}, \ \ c \bl(V)=0,\ \ d(v\otimes t^r)=rv\otimes t^r.\end{gather*} \noindent Clearly $\bl$ takes finite--dimensional $L(\lie g)$--modules to integrable\index{module!integrable} $\hat L(\lie g)$--modules. It will become clear from the results of the next section, that the functor $\bl$ in general sends irreducible finite--dimensional  modules to completely reducible integrable\index{module!integrable} ones. In  \cite{CG0} it is shown that $\bl$  preserves enough of the homological properties of $L(\lie g)$--modules to make it worthwhile to restrict one's attention to modules for $L(\lie g)$.

\subsection{}\label{C:section 4.4} We digress, briefly, from our study of $\cal I$ to define the imaginary integral root vectors and the integral form of $\bu(\hat L(\lie g))$. These were introduced by Garland in \cite{Ga} and allow us to study lattices in integrable\index{module!integrable}  modules and hence also the representations of the hyperalgebra\index{hyperalgebra} of $\hat L(\lie g)$. It is also worthwhile to note that the integral imaginary root vectors  play a major role in the vertex operator construction \cite{FJ} of the basic representation\index{basic representation} which was defined in Section \ref{C:section 3.8}.

As in the case of $\lie g$ we let $\bu_{\bz}(\hat L(\lie g))$ be the $\bz$--subalgebra of $\bu(\hat L(\lie g))$ generated by the elements $(x_\alpha^\pm\otimes t^s)^{(r)}$, where $\alpha\in \Phi^+$, $s\in\bz$ and $r\in \bn$. We know from Section \ref{C:section 2.11} that one has to be more careful with  the  Cartan type elements $(h_i\otimes t^s)$, $s\in\bz$. For $s=0$, $r\in\bn$ and $1\le i\le n$,  the elements $\binom{h_i}{r}$ are in $\bu_\bz(\lie g)$ and hence in $\bu_{\bz}(\hat L(\lie g))$.

To define the analogs of the divided powers for the imaginary root vectors  $h_i\otimes t^s$, $i\in I$, $s\in\bz$, $s\ne 0$ one tries to modify suitably  the principle we used in Section \ref{C:section 2.11} (see equation \eqref{C:zform}). Consider the  triangular decomposition $$\bu(\hat L(\lie g))=\bu(\hat{\lie n}^-)\bu(\hat{\lie h})\bu(\hat{\lie n}^+),$$ and note that for all $\alpha\in \Phi^+$, $r,m\in \bn$ the element $(x_\alpha^+\otimes t)^{(r)}(x^-_\alpha)^{(m)}$, is not in the order prescribed on the right hand side. In \cite{Ga}, Garland proved a remarkable formula which rewrites this element in the correct order. We will not reproduce his entire  formula here, but content ourselves with writing down the correct analog (which emerges from his formula)  of the divided power of the imaginary root vectors.
For each $i\in I$ we define   formal power series $P^\pm_i(u)$ in an indeterminate $u$ with values in the commutative algebra $\bu(\lie h\otimes t^{\pm 1}\bc[t^{\pm 1}])$ by \begin{equation}\label{C:intimgvec} P^\pm_i(u)=\sum_{s=0}^\infty P^\pm_{i,s}u^s={\rm{exp}}\left(-\sum_{k=1}^\infty \frac{h_i\otimes t^{\pm k}}{k}\right),\end{equation} and we can now state Garland's theorem.
\begin{thm} Fix an order on the set  $$\{(x_\alpha^\pm\otimes t^s)^{(r)}:\alpha\in \Phi^+, s\in\bz, r\in\bn\}\cup\{P^\pm_{i,s}: 1\le i\le n, s\in\bn\}\cup\{\binom{h_i}{r}: i\in I, r\in\bn\}. $$  The ordered monomials from this set are a $\bz$--basis of $\bu_{\bz}(\hat L(\lie g))$ and a $\bc$--basis of $\bu(\hat L(\lie g))$. Moreover, if $\lambda\in\hat P^+$ and $\hat v_\lambda\in\hat V(\lambda)_\lambda$, then $\bu_{\bz}(\hat L(\lie g))\hat v_\lambda$ is a $\bz$--lattice in  $\hat V(\lambda)$.\hfill\qedsymbol
\end{thm}

The twisted version of this theorem   may be found in \cite{Mitz}.

\subsection{}\label{C:section 4.5} In view of this theorem it would be natural to pass to a field $\Bbb F$ of characteristic $p$ as we did in Section \ref{C:charpg} and study the representations  $\bu(\hat L(\lie g))_\Bbb F$ and in particular the modules $\hat V(\lambda)_\Bbb F$. However, this direction does not seem to have been pursued and except for the result we now discuss,  virtually nothing is known about the irreducible modules of positive level in characteristic $p$. In \cite{CJ}, a quantum analog of this problem was studied and it was shown there that the basic representation\index{basic representation} of the quantum affine\index{affine Lie algebra!quantum} algebra remained irreducible at an $m^{th}$  root of unity, where $m$ is odd and coprime to the determinant of the Cartan matrix of $\lie g$. This was proved essentially by constructing the basic representation\index{basic representation} at a root of unity  and showing that it had the same character. The proof given in that paper works verbatim in the characteristic $p$ case as long as $p$ is odd and coprime to the determinant of the Cartan matrix of $\lie g$.
The case of level zero modules in characteristic $p$ has been studied by Jakelic and Moura in \cite{JM1}, \cite{JM2}, \cite{JM3} and we shall discuss this in  the next sections.

      \section{Finite--dimensional modules for loop algebras and their\break \ generalizations}\label{C:section 5}
In this section, we discuss the category of finite--dimensional representations  of loop algebras, study extensions between irreducible modules and describe the blocks\index{block} of this category.

\subsection{}\label{C:section 5.1} For $a\in\bc^\times$ let $\ev_a: L(\lie g)\to\lie g$ be the homomorphism of Lie algebras given by $$\ev_a(x\otimes f)=f(a)x.$$ If $V$ is a $\lie g$--module, we denote the corresponding $L(\lie g)$-module  by $\ev_a^*V$\index{module!evaluation}. Together with the discussion in Section \ref{C:section 4.3}, this explains the notation  $L(\ev_a^*V)$ used in Section \ref{C:section 4.2}. The following result classifies irreducible finite--dimensional representations of $L(\lie g)$.
\begin{thm} \label{C:findirr}An irreducible finite--dimensional  representation of $L(\lie g)$ is isomorphic to a tensor product $\ev^*_{a_1}V(\lambda_1)\otimes\cdots\otimes \ev^*_{a_k} V(\lambda_k)$, where $\lambda_s\in P^+$, $a_s\in\bc^\times$ with $a_s\ne a_r$ for $1\le s, r\le k$. Conversely any such tensor product is irreducible.\hfill\qedsymbol
\end{thm}
This result can be deduced easily from Proposition \ref{C:loopirr} and Theorem \ref{C:inv}(iii). A more direct proof was given by Rao in \cite{Rao1} along the following lines. If $\rho:L(\lie g)\to\End V$ is a finite--dimensional representation, then the kernel of $\rho$ is an ideal of finite--codimension in $L(\lie g)$. It is not too hard to show that any ideal in $L(\lie g)$ must be of the form $\lie g\otimes I(V)$ where $I(V)$ is an ideal in $\bc[t,t^{-1}]$. If $V$ is irreducible then one proves that $I(V)$ is a product of distinct maximal ideals and hence is generated by an element $(t-a_1)\cdots (t-a_k)$ for some $k\in\bn$ and distinct elements $a_1,\cdots, a_k$ in $\bc^\times$
Moreover $\bc[t,t^{-1}]/I(V)$ is a vector space of dimension $k$  and it follows now that $L(\lie g)/(\lie g\otimes I(V))$ is isomorphic to a direct sum of $k$--copies of $\lie g$ and hence is a semisimple Lie algebra. The theorem now follows by using the representation theory of semisimple Lie algebras.

 For the twisted affine Lie algebras the corresponding result was proved in \cite{CFS}. To describe it, note that since $L(\lie g,\sigma, m)$ is a subalgebra of $L(\lie g)$, one can regard the tensor product $V= \ev^*_{a_1}V(\lambda_1)\otimes\cdots\otimes \ev^*_{a_k} V(\lambda_k)$ as a module for the twisted algebra. In general it is not irreducible as a module for the twisted algebra, but if one imposes the additional condition that $a_s^m\ne a_r^m$ for all $1\le s,r\le k$, where $m$ is the order of $\sigma$, then $V$ is an irreducible module for $L(\lie g,\sigma, m)$. Moreover, one can also prove that these are exactly (up to isomorphism) all the finite--dimensional irreducible modules.

\subsection{}\label{C:section 5.2} An element  $a\in \bc^\times$ determines a maximal ideal $(t-a)$ of $\bc[t,t^{-1}]$ and another way to formulate  Theorem \ref{C:findirr} is to say that the irreducible finite--dimensional representations are parametrized by finitely--supported functions from the maximal spectrum\index{maximal spectrum} of $\bc[t,t^{-1}]$ to $P^+$, where by finitely supported, we mean $\chi(M)=0$ for all but finitely many maximal ideals.  We  generalize the theorem to the case of $\lie g\otimes A$, where $A$ is a finitely generated commutative associative algebra algebra  over $\bc$. The space $\lie g\otimes A$ is a Lie algebra as usual: $[x\otimes a, y\otimes b]=[x,y]\otimes ab$.  Let $\max A$ be the set of maximal ideals in $A$ and denote by  $\Xi(\max A,P^+)$ the set of  finitely supported functions from $\max A$ to $P^+$,  and set$$\supp\psi=\{M\in\max A:\psi(M)\ne 0\},\ \ \psi\in\Xi(\max A, P^+).$$ For $M\in\max A$, let $\ev_M:\lie g\otimes A\to \lie g$ be the map of Lie algebras which is  induced by the algebra homomorphism $A\to A/M\cong\bc$ and we denote by $\ev_M^*V(\lambda)$ the representation of $\lie g\otimes A$ obtained by pulling back $V(\lambda)$ through $\ev_M$.
\begin{thm}\label{C: findimgen} Let $A$ be a finitely generated commutative and associative algebra over $\bc$. For  $\chi\in\Xi(\max A, P^+)$, the $\lie g\otimes A$--module  $\bigotimes_{M\in\supp\chi}\ev_M^*V(\chi(M))$ is irreducible and conversely any finite dimensional irreducible representation is isomorphic  to one of these.\hfill\qedsymbol\end{thm}
The theorem has been proved by many people in various cases: in the case when $A$ is the polynomial ring in $r$ variables the result was first proved by Rao \cite{Rao1}, the general case given above appears in \cite{CFK}. In \cite{Lau}, Lau also obtains a proof of this statement for the Laurent polynomial ring in $r$ variables and he is also able prove the result in the  twisted case. It is important to note here, that unlike in the case when $A=\bc[t,t^{-1}]$, the twisted algebras are not determined by a Dynkin diagram automorphism and  one has many non--isomorphic twisted algebras and Lau works in this generality. The methods of all these papers are  algebraic and similar to the proof sketched above for the loop algebra. A more geometric approach is developed in  \cite{NSS} and their methods allow them to develop a uniform approach which works for twisted and untwisted algebras of the form $\lie a\otimes A$ where $\lie a$ is an arbitrary Lie algebra. In this generality  the irreducible representations are not always tensor products of evaluations, but they are very nearly so, up to tensoring with a one--dimensional representation.

\subsection{}\label{C:section 5.3} An obvious and not very difficult question that one should ask is what kind of module one gets if one allows  $a_r=a_s$ for some pair $1\le r, s,\le k$ in the statement of Theorem \ref{C:findirr}. The following proposition proved in \cite{CGfree} shows that the module is completely reducible and more generally gives a necessary and sufficient condition    for a finite--dimensional module to be completely reducible. Recall that for any $L(\lie g)$--module $V$ we have defined in Section \ref{C:findirr} an ideal $I(V)$ of $\bc[t,t^{-1}]$ which is maximal with the property that $(\lie g\otimes I(V))V=0$. More precisely, $\lie g\otimes I(V)$ is the annihilating ideal of $V$  in $L(\lie g)$.
\begin{prop} Let $\lambda,\mu\in P^+$ and $a\in\bc^\times$. Then $\ev^*_a V(\lambda)\otimes\ev^*_a V(\mu)$ is isomorphic to a direct sum of submodules of the form $\ev^*_a V(\nu)$, $\nu\in P^+$  and $$\Hom_{L(\lie g)}(\ev^*_a V(\nu), \ev^*_a V(\lambda)\otimes\ev^*_a V(\mu))\cong\Hom_{\lie g}(V(\nu), V(\lambda)\otimes V(\mu)).$$ More generally, a finite--dimensional module  $V$ is completely reducible iff the ideal $I(V)$ is a product of distinct maximal ideals. Analogous statements hold for algebras of type $\lie g\otimes A$.
\end{prop}
\subsection{}\label{C:section 5.4} We now describe extensions between irreducible $L(\lie g)$--modules. We begin with a simple case and compute $\Ext_{L(\lie g)}^1(\ev^*_a V(\lambda), \ev^*_a V(\mu))$ for $a\in\bc^\times$.

Let \begin{equation}\label{C:se1} 0\to\ev^*_a V(\mu)\to U\to \ev_a^* V(\lambda)\to 0\end{equation} be  a  short exact sequence of $L(\lie g)$--modules.  Regarded as a short exact sequence of $\lie g$--modules this sequence is split and we can pick $u_\lambda\in U_\lambda$ with $\bu(\lie g)u\cong V(\lambda)$.  If \eqref{C:se1} is not split as a short exact sequence of $L(\lie g)$--modules,  there  must exist $r\in\bz$ such that the   $\lie g$--module map,$$\lie g\otimes \bu(\lie g)u_\lambda\to U, \qquad (xt^r\otimes gu_\lambda)\to xt^rgu_\lambda,$$ has a non--zero projection (as a $\lie g$--module) onto the image of $\ev^*_a V(\mu)$, i.e., $$\Ext_{L(\lie g)}^1(\ev^*_a V(\lambda), \ev^*_a V(\mu))\ne 0\implies \Hom_{\lie g}(\lie g\otimes V(\lambda), V(\mu))\ne 0.$$
Conversely, given $\pi:\lie g\otimes V(\lambda)\to V(\mu)$ a map of $\lie g$--modules it is proved in \cite{CM} that  the formula
$$xt^r(v,w) = (a^rxv, a^rxw + ra^{r-1}\pi(x
 v)),\ \ v\in V(\lambda), w\in V(\mu),\ \ a\in\bc^\times,\ \ r\in\bz, \ \ x\in\lie g,$$ defines an indecomposable\index{indecomposable} $L(\lie g)$--module denoted $V(\lambda,\mu,a)$ on $V(\lambda)\oplus V(\mu)$ and hence one can conclude that $$\Ext_{L(\lie g)}^1(\ev^*_a V(\lambda), \ev^*_a V(\mu))\cong \Hom_{\lie g}(\lie g\otimes V(\lambda), V(\mu)),$$ which can be shown to be equivalent to $$\Ext_{L(\lie g)}^1(\ev^*_a V(\lambda), \ev^*_a V(\mu))\cong \Hom_{L(\lie g)}(\ev^*_a\lie g\otimes \ev^*_a V(\lambda), \ \  \ev^*_a V(\mu)).$$
 The general statement proved in \cite{CGfree} (for $\lie g\otimes \bc[t]$ but which has an obvious modification to $L(\lie g)$) is:
 \begin{thm}\label{C:extirr} Let $V$ and $V'$ be irreducible finite--dimensional $L(\lie g)$--modules. Then $$\dim\Ext^1_{L(\lie g)}
(V, V ') =
\sum_
{a\in\bc^\times}\dim\Hom_{L(\lie g)}(\ev_a^*\lie g\otimes V,
 V').$$\hfill\qedsymbol
\end{thm}
This condition can be made quite explicit using the results stated above, namely if we write $V$ and $V'$ as a tensor product of evaluation representations as in Theorem \ref{C:findirr} one can compute the dimension of $\Ext_{L(\lie g)}^1(V, V')$ as the sum of multiplicities  of the adjoint representation of $\lie g$ in a tensor product of the form $V(\lambda)\otimes V(\mu)^*$ for suitable $\lambda$, $\mu$.

\subsection{}\label{C:section 5.5} The analogous result has not been established as yet for the twisted Lie algebras. It has been solved in \cite{Ko} for the more general Lie algebra $\lie g\otimes A$ and in fact the main theorem in this paper  is written in a very  explicit way. There is one new feature though which is not seen in the case of $\bc[t,t^{-1}]$ which we explain again in a simple case. For $\lambda,\mu\in P^+$ and $M\in\max A$, we have $$\Ext^1_{\lie g\otimes A}(\ev_M^*V(\lambda), \ev^*_M V(\mu))\cong\Hom_{\lie g}(\lie g\otimes V(\lambda), V(\mu))\otimes {\rm{Der}}(A, A/M),$$ where ${\rm{Der}}(A, A/M)$ is the space of all linear maps $A\to A/M$ satisfying $D(ab)= aD(b)+bD(a)$. In the case when $A=\bc[t,t^{-1}]$ this space is one--dimensional and so is invisible.

\subsection{}\label{C:section 5.6} To continue our study of finite--dimensional representations, we look back at Section 2 and notice  that one of the natural things to do is to understand the blocks\index{block} of the category. As we already noted in Section 4, the theory of central characters is not available for affine or loop algebras. However, it is still possible to describe the blocks\index{block} and  to do this one uses a family of universal indecomposable\index{indecomposable} finite--dimensional modules introduced in \cite{CPweyl} and called the Weyl modules\index{module!Weyl} for affine Lie algebras. These are not the same modules as the ones discussed in Section 2, but are in some sense obtained in the same way. They can be regarded as the classical limit of irreducible representations of quantum affine algebras and we shall say more about this later. The first step is to organize the irreducible representations in some nice way and for this, we  recall from \cite{CM} the notion of the  spectral character of a $L(\lie g)$--module.
\begin{defn} Given an irreducible representation $V=\otimes_{s=1}^k \ev^*_{a_s}V(\lambda_s)$ of $L(\lie g)$ the spectral character of $V$ is defined to be the function $\chi_V:\bc^\times\to P/Q$ given by$$\chi_V(z)=\begin{cases} 0,\ \ z\notin\{a_1,\cdots,a_k\},\\ \lambda_s+Q,\ \ z=a_s, \ \ 1\le s\le k.\end{cases}$$ A finite dimensional $L(\lie g)$--module is said to have spectral character if all the  irreducible modules occurring in a Jordan--H\"older series have the same spectral character.\hfill\qedsymbol
\end{defn}
\noindent Note that it is very easy for two modules to have the same spectral character, for instance if $\lambda,\mu\in P^+\cap Q$, then $$\chi_{\ev^*_aV(\lambda)}=\chi_{\ev^*_bV(\mu)},\ \  {\rm{for\ all}}\ \ a,b\in\bc^\times.$$

\subsection{}\label{C:section 5.7}\label{C:blockdecomp} We recall the definition of the blocks of a category $\cal C$\index{block}. \begin{defn} Say that two indecomposable\index{indecomposable} objects $U,V\in\Ob\cal C$ are linked\index{linked} if there do not exist abelian
subcategories $\cal C_k$, $k =1,2$ such that $\cal C = \cal C_1\oplus\cal C_2$ with $U\in\cal C_1$ and $V\in\cal C_2$. If $U$ and $V$ are decomposable then
we say that they are linked if every indecomposable\index{indecomposable} summand of $U$ is linked to every indecomposable\index{indecomposable}
summand of $V$ .
This defines an equivalence relation on $\cal C$ and the equivalence classes are called the  blocks\index{block} of $\cal C$.
\end{defn}

The main result proved in \cite{CM} is the following theorem which describes the blocks\index{block} of the category of finite--dimensional representations.
\begin{thm} Let $V$ be a finite--dimensional $L(\lie g)$ module. Then $V$ is isomorphic to a direct sum of finite--dimensional  modules $V_r$, $1\le r\le k$ each of which admits a spectral character. Moreover, any two modules with the same spectral character are linked.
\end{thm}
As a consequence of this theorem, we see that the blocks\index{block} are very large subcategories, for instance when $\lie g$ is of type $E_8$ there is only one block\index{block}, i.e., any two finite--dimensional modules are linked.
A similar result was proved in \cite{S} for the twisted algebras $L(\lie g,\sigma, m)$ and in \cite{Ko} for the Lie algebras $\lie g\otimes A$. All the proofs are similar and involve two ingredients. The first is the module $V(\lambda,\mu,a)$ constructed in Section \ref{C:section 5.4} and the second are the Weyl modules\index{module!Weyl} which we  discuss in the next section.

\subsection{}\label{C:section 5.8} Most of the results of this section  have analogs in positive characteristic  for the   hyperalgebras\index{hyperalgebra} $\bu(L(\lie g))_{\Bbb F}$ and these are studied in \cite{JM1}, \cite{JM2}. The one exception is the result on extensions  between irreducible modules. This is because the proof given in characteristic zero relies on the fact that a finite--dimensional representation of a complex simple Lie algebras is completely reducible, which as we remarked in Section \ref{C:charpg} is false in positive characteristic.
\section{Weyl modules, restricted Kirillov--Reshetikhin and beyond}\label{C:section 6}
We begin this section by discussing briefly the input coming from the representation theory of quantum affine algebras, which has helped to identify interesting families of finite--dimensional representations of loop algebras. The interested reader could consult \cite{CH} and the references therein for further details on representations of quantum affine algebras.
\subsection{}\label{C:section 6.1} In the 1980's Drinfeld and Jimbo introduced the quantized enveloping algebras $\bu_q(\lie a)$ of a Kac--Moody Lie algebra $\lie a$. These algebras depend on a parameter $q$ which can be either a complex number, in which case one regards $\bu_q(\lie a)$ as an algebra over $\bc$ or a formal variable in which case one regards it as an algebra over the function field $\bc(q)$, and informally speaking if we put $q=1$ we get the usual enveloping algebra of $\lie a$. In the case when $\lie g$ is simple, Lusztig \cite{Lus88} proved that for generic $q$, the irreducible finite--dimensional representations of $\bu_q(\lie g)$ are given by elements  of $P^+$ and moreover that the module associated to $\lambda\in P^+$ has the same character (suitably understood) as the module in the $q=1$ case, namely it is the character of  the $\lie g$--module $V(\lambda)$. An analogous statement was also proved  for the positive level integrable\index{module!integrable} modules of the quantized affine algebras. In the case when $q$ is a primitive $m^{th}$--root of unity, the picture resembles the case of positive characteristic and the irreducible modules are in general smaller.

\subsection{}\label{C:section 6.2} In the case of finite--dimensional modules for the quantized affine algebras the difference between generic $q$ and the roots of unity case is already visible when $q=1$. In other words the action of $L(\lie g)$ on evaluation representations and their tensor products given in Theorem \ref{C:findirr} cannot be deformed to give an action of the quantum affine algebra on the same space. Or, to put it in yet another way, the $q=1$ limit of irreducible representations of quantum affine algebras generally give rise to a reducible indecomposable\index{indecomposable} representation of $L(\lie g)$. This was the motivation for introducing in \cite{CPweyl} the definition of Weyl modules for affine Lie algebras. We should remark here that this phenomenon was first noted in \cite{Dr} in connection with the representation theory of the closely related Yangians, which are deformations of the current algebra\index{current algebra}, $\lie g\otimes\bc[t]$.

\subsection{}\label{C:section 6.3} The classification of irreducible finite--dimensional representations for quantum affine algebras was obtained in \cite{CPqa} and is similar to the one obtained by Drinfeld for the Yangians. It was also very similar to the abstract classification (rather than the more explicit one described in these notes) of loop modules given in \cite{Cinv}. However, the fact that there is no analog of the evaluation homomorphism from the quantum loop algebra to the quantum simple algebra has meant that describing the irreducible modules in a concrete way has been challenging and a variety of methods have been developed to understand these modules better. One of these methods is to understand the multiplicity with which an irreducible module for the quantum simple algebra occurs in an irreducible module for the quantum loop algebra. This can be turned into a problem of the non--quantum case for the following reason. It was proved in \cite{CPweyl} that under natural conditions the irreducible representations of the quantum affine algebra admitted a $\bc[t,t^{-1}]$--form and hence could be specialized to $q=1$ and gives a representation of $L(\lie g)$. A result of Lusztig implies that understanding the $\lie g$--module decomposition of this module is the same as answering the question for the quantum case. With this motivation out of the way, we now return our focus to representations of $L(\lie g)$.

\subsection{}\label{C:section 6.4} For $\lambda\in P^+$, let $W(\lambda)$ be the $L(\lie g)$--module generated by an element $w_\lambda$ with relations $$h w_\lambda=\lambda(h)w_\lambda,\ \ L(\lie n^+) w_\lambda=0,\ \ (x^-_{\alpha_i})^{\lambda(h_{\alpha_i})+1}w_\lambda=0.$$ We call $W(\lambda)$ the global Weyl module\index{module!Weyl!global}. The following is elementary.
 \begin{lem} The module $W(\lambda)$ is integrable\index{module!integrable}  for all $\lambda\in P^+$. Moreover  if $V$ is any integrable module, then $$\Hom_{L(\lie g)}(W(\lambda), V)\cong V_\lambda^+,\qquad \ V_\lambda^+=\{v\in V_\lambda: L(\lie n^+)v=0\}.$$
  In particular, if $\ev^*_{a_1}V(\lambda_1)\otimes\cdots\otimes \ev^*_{a_k}V(\lambda_k)$ is an irreducible finite--dimensional $L(\lie g)$--module with $\sum_{s=1}^k \lambda_s=\lambda$, then $$\dim\Hom_{L(\lie g)}(W(\lambda),\ev^*_{a_1}V(\lambda_1)\otimes\cdots\otimes \ev^*_{a_k}V(\lambda_k))=1.$$

\hfill\qedsymbol
 \end{lem}
 \noindent The lemma shows immediately that the  module $W(\lambda)$ is infinite--dimensional if $\lambda\ne 0$ and in fact that $W(\lambda)_\lambda$ is infinite--dimensional.  For otherwise, the action of the algebra $L(\lie h)$ on $W(\lambda)_\lambda$ would be a direct sum of a finite number of generalized eigenspaces. On the other hand the lemma says that $W(\lambda)$ has infinitely many irreducible quotients $\ev^*_{a_1}V(\lambda_1)\otimes\cdots\otimes \ev^*_{a_k}V(\lambda_k)$ and a simple calculation shows that the eigenvalues of $L(\lie h)$ on the image of $w_\lambda$ in each of these quotients is different.

 \subsection{}\label{C:section 6.5} Let us look at a further consequence of this lemma. Recall the modules $V(\lambda,\mu,a)$ defined in Section \ref{C:section 5.4}. If $\lambda-\mu\in Q^+$, then with a little computation, one can prove that $V(\lambda,\mu,a)$ is generated as a $L(\lie g)$--module by the element $v_\lambda\in V(\lambda)$ and that $L(\lie n^+)v_\lambda=0$ and $\dim V(\lambda,\mu,a)_\lambda=1$, if $\lambda\ne \mu$. Hence $\Hom_{L(\lie g)}(W(\lambda), V(\lambda,\mu,a))\ne 0$ and so $W(\lambda)$ also has cyclic reducible (indecomposable\index{indecomposable}) finite--dimensional quotients. It is reasonable therefore,  to ask if there are  a family of maximal or universal  finite dimensional quotients of $W(\lambda)$: a family of finite--dimensional quotients $V$ which are cyclic and have $\dim V_\lambda=1$ and are  such that any other such finite--dimensional quotient of $W(\lambda)$ is a quotient of one of these. The answer is yes and these are called the local Weyl modules\index{module!Weyl!local} and we shall spend some time discussing these modules. We remark that it is these local Weyl modules which appear as the  $q=1$ limits of representations of quantum affine algebras.

\subsection{}\label{C:section 6.6} The module $W(\lambda)$ admits a right $L(\lie h)$--module structure, given by $$(yw_\lambda)u= yuw_\lambda,\ \ y\in\bu(L(\lie g)),\ \ u\in L(\lie h),$$ and hence $W(\lambda)$ is a $(L(\lie g), L(\lie h))$--bimodule.  Given $\chi:L(\lie h)\to\bc$, set $$W(\lambda,\chi)=W(\lambda)\otimes_{L(\lie h)} \bc_\chi,$$ where $\bc_\chi$ is the one--dimensional $L(\lie h)$--module defined by $\chi$.
The following was proved in \cite{CPweyl}.
\begin{thm} Let  $\lambda\in P^+$ and $\chi\in L(\lie h)^*$. Then the module $W(\lambda,\chi)$ is finite--dimensional and is non--zero iff there exists $k\in\bn$, $\lambda_1,\cdots,\lambda_k\in P^+$ and $a_1,\cdots,a_k\in\bc^\times$ such that $$\chi(h\otimes t^s)=\sum_{j=1}^k\lambda_j(h)a_j^s.$$ Moreover if $(\lambda,\chi)$ and $(\lambda',\chi')$ are such that $W(\lambda,\chi)$ and $W(\lambda',\chi')$ are non--zero, then $$W(\lambda,\chi)\cong W(\lambda',\chi')\ \ \iff \ \ \lambda=\lambda',\ \ \chi=\chi'.$$\hfill\qedsymbol
\end{thm}
The notation used here, although different from that used in \cite{CPweyl}, is consistent with the notation of the previous sections of these notes.

\subsection{}\label{C:section 6.7} The local Weyl modules\index{module!Weyl!local} have a tensor product factorization proved in \cite{CPweyl} which is  similar to that given in Theorem \ref{C:findirr}. For $\lambda\in P^+$ and $a\in\bc^\times$, let $\chi_{\lambda,a}\in L(\lie h)^*$ be given by $\chi_{\lambda,a}(h\otimes t^r)=a^r\lambda(h)$.
\begin{prop} Any non--zero local Weyl module is isomorphic to a tensor product  of the form $W(\lambda_1, \chi_{\lambda_1,a_1})\otimes\cdots\otimes W(\lambda_k, \chi_{\lambda_k,a_k})$ for some $\lambda_s\in P^+$, $a_s\in\bc^\times$, $1\le s\le k$ with $a_j\ne a_m$ if $j\ne m$.\hfill\qedsymbol
\end{prop}
As a result of the proposition, to understand the $\lie g$--module structure or the dimension of the local Weyl modules it suffices to study the modules $W(\lambda,\chi_{\lambda,a})$ for $(\lambda,a)\in\bc^\times$. These questions are addressed in \cite{CPweyl},\cite{CL},\cite{FoL} for $\lie g$ of type $A_1$, $A_n$ and $A,D,E$ respectively and in these cases the local  Weyl modules can be identified with Demazure\index{module!Demazure} modules in positive level representations of affine Lie algebras. This is definitely false for non--simply laced algebras where the local Weyl modules are bigger than the Demazure\index{module!Demazure} modules. In general the dimension of the local Weyl\index{module!Weyl!local} modules can be deduced from the work of \cite{BN}. As a consequence one knows that the dimension of $W(\lambda,\chi_{\lambda,a})$ is independent of the choice of $a\in\bc^\times$ which in turn can be used to prove that $W(\lambda)$ is free as a right module for a suitably defined polynomial ring (the quotient of $\bu(L(\lie h))$ by the annihilator of $w_\lambda$).

\subsection{}\label{C:section 6.8} The Weyl modules\index{module!Weyl!local} for $\lie g\otimes A$ have been defined and studied in \cite{CFK} for an arbitrary associative, commutative algebra $A$. The results of that paper  show that there are many important differences between the study of Weyl modules for the loop algebra and the more general case. As in the case for irreducible modules, one works with maximal ideals, and the local Weyl modules for a fixed $\lambda$ are parametrized by a maximal ideal of $A$.  The dimension of the local Weyl modules even in the simplest cases does depend on the maximal ideal  even when $A$ is the polynomial ring in two variables and hence these modules are no longer free right modules. There are thus many interesting problems that arise from passing to the higher dimensional cases.

\subsection{}\label{C:section 6.9} Local Weyl\index{module!Weyl!local} modules have also been studied in positive characteristic in \cite{JM1},\cite{JM2}. As in the case of simple Lie algebras, the definition needs modification, the most crucial one being that  the conditions on $L(\lie h)$ must be replaced with the conditions on  the integral imaginary root vectors defined by Garland. The authors also study the question of base change, i.e., the functor of extension of scalars from the category of finite-dimensional $\bu(L(\lie g))_{\Bbb K}$-modules to that of $\bu(L(\lie g))_{\Bbb F}$-modules. For simple Lie algebras, the functor induces a bijection from the set of isomorphism classes of highest-weight\index{module!highest weight} modules for $\bu(\lie g)_\mathbb K$ to that of isomorphism classes of highest-weight modules for $\bu(\lie g)_\mathbb F$, where $\mathbb F$ is an algebraic extension of $\mathbb K$.  In particular, the characters\index{character} of the irreducible modules and the Weyl modules\index{module!Weyl} for the hyperalgebras\index{hyperalgebra} $\bu(\lie g)_\Bbb L$ do not depend on the choice of the field $\mathbb L$, but only on its characteristic.
In the case of hyper loop algebras, the
story is more complicated since the functor does not send irreducible modules to
irreducible
modules and, hence, does not preserve the length of a module. For more details on this functor in this case see \cite{JM3}.

\subsection{}\label{C:section 6.10} An important family of quotients of the Weyl modules are the Kirillov--Reshetikhin (KR)\index{module!KR}--modules. In the quantum case, these modules are very important, since they admit a crystal basis. Most finite--dimensional representations do not admit crystal bases. These modules first arose from the study of solvable lattice models and a good reference is \cite{HKOTY} and again more recent references can be found in \cite{CH}. The results of \cite{CPweyl} show that one can put $q=1$ and get modules for the loop algebra. The KR--modules, for our purposes, can be viewed as quotients of the local Weyl\index{module!Weyl!local} modules $W(m\omega_i, \chi_{m\omega_i, a})$, where $i\in I$, $m\in\bn$ and $a\in\bc^\times$. They are obtained by imposing a single additional relation: $$(x_{\alpha_i}^-\otimes (t-a)w_{m\omega_i})=0.$$ Amongst other things, one of the interesting problems is to compute the $\lie g$--character of this module and this is done in \cite{CM} for the classical Lie algebras and for some nodes of the exceptional Lie algebras. These modules come equipped with a $\bn$--grading and, in fact, the graded character is computed in \cite{CM} and shown to coincide with the conjectural character formulae in \cite{HKOTY}. The relationship between the grading and the parameter $q$ that appears in \cite{HKOTY}, though, is quite mysterious. An obvious question from a mathematical point of view, is what happens when $m\omega_i$ is replaced with an arbitrary weight\index{weight} $\lambda\in P^+$. Partial answers to this question can be found in \cite{M}. An alternative way to think of KR-modules\index{module!KR} is as projective\index{module!projective} modules in a suitable subcategory and this is being explored in ongoing work with Jacob Greenstein.

\section{Koszul algebras, quivers, and highest weight categories}\label{C:section 7}

In this section we discuss an approach, developed jointly with Jacob Greenstein in  \cite{CG1}, \cite{CG2}, to the graded  representation theory of the algebras $\lie g\otimes\bc[t]$. We  follow the axiomatic approach of Cline, Parshall and Scott which  has been successful in both the study of $\cal O$ and representations in characteristic $p$. The basic idea is to look at subcategories, each of which has finitely many simple objects and enough projectives (or injectives), and to relate it to the module category of a finite--dimensional algebra. One can then also study the quiver with relations attached to such a subcategory. It is usual to assume that the subcategory is a block\index{block} and also that one has a natural ordering of the simple objects in the block. This is usually given by the Bruhat order on the Weyl group\index{Weyl group}, which in turn is related  to extensions between Verma modules.

 However, in our case, as we have seen in Section \ref{C:blockdecomp}, the blocks\index{block} of the finite--dimensional representations are very large and this remains true if we restrict our attention to graded representations of $\lie g\otimes\bc[t]$.
 Nevertheless, we are able to identify interesting subcategories, which contain finitely many  irreducibles and the corresponding (generalized) Kirillov--Reshetikhin modules discussed in  Section 6. We are also able to define an order on the simple modules which comes  from an  understanding of the extensions between simple modules.   The endomorphism algebra of such a subcategory is a finite--dimensional algebra given by a quiver with relations. Under further restrictions  the endomorphism algebra  has a natural grading and using the results  of \cite{BGS} we prove that the grading is Koszul. In many cases (in fact for an infinite family) it is possible \cite{G} to compute the relations explicitly and we give some concrete examples below. The twisted analogs of these results and in fact further generalizations may be found in \cite{CKR}.

\subsection{}\label{C:section 7.1} The algebra $\lie g[t]=\lie g\otimes\bc[t]$ is naturally graded by $\bn$ with the $r^{th}$--graded piece being $\lie g\otimes t^r$ for $r\in\bn$. The enveloping algebra $\bu(\lie g[t])$ acquires a well--defined natural grading by $\bn$ as well: we say that an element of $\bu(\lie g[t])$ has grade $s$ if it is a linear combination of elements of the form $(x_1\otimes t^{s_1})\cdots (x_k\otimes t^{s_k})$ where  $k\in\bn$, $s_j\in\bn$, $x_j\in\lie g$, $1\le j\le k$ and $s=s_1+\cdots +s_k$. Let $\lie g[t]_+=\lie g\tensor t\bc[t]$, which is clearly a Lie ideal of $\lie g[t]$. Note that $\bu(\lie g[t])\cong
\bu(\lie g[t]_+)\tensor \bu(\lie g)$ and that the graded pieces $\bu(\lie g[t]_+)[r]$ are finite dimensional $\lie g$-modules.
Hence for all $\lambda\in P^+$, we see that $\bu(\lie g[t])[s]\tensor_{\bu(\lie g)} V(\lambda)$ is a finite--dimensional $\lie g$--module.

 Let $\cal G$ be the category
 whose objects are $\bn$--graded modules $V$ with finite dimensional graded pieces,
 \begin{gather*} V=\bigoplus_{r\in\bn} V[r],\ \ \text{all}\; \dim V[r]<\infty,\\ \ (\lie g\otimes t^s)V[r]\subset V[r+s],\  \ r,s\in\bn,\end{gather*} and the morphisms are given by $$\Hom_{\cal G}(V,W)=\{f\in\Hom_{\lie g[t]}(V,W): f(V[r])\subset W[r],\ \ r\in\bn\}.$$  For $r\in\bn$ let $\tau_r$ be the grading shift functor
 $\cal G\to \cal G$, that is for all
 $V\in\Ob\cal G$, $\tau_r V$ is isomorphic to $V$ as a $\lie g[t]$--module and  $(\tau_rV)[s]=V[r+s]$.

\subsection{}\label{C:section 7.2} Let  $\ev_0:\lie g[t]\to \lie g$, given by $x\otimes t^r\to\delta_{r,0}x$, be the evaluation at zero.
For $\lambda\in P^+$ we can regard the $\lie g[t]$--module $\ev^*_0V(\lambda)$ as an object of $\cal G$ by considering it as being concentrated in grade zero.
Set $V(\lambda,r)=\tau_r \ev^*_0 V(\lambda)$.
The projective cover of $V(\lambda,r)$ in $\cal G$ is $$P(\lambda,r)=\tau_r\bu(\lie g[t])\otimes_{\bu(\lie g)}V(\lambda),\qquad P(\lambda,r)[s]=\bu(\lie g[t])[s-r]\otimes_{\bu(\lie g)}V(\lambda),\qquad s\ge r.$$
The following is not hard to prove.
\begin{prop} Any irreducible object in $\cal G$ is isomorphic  to $V(\lambda,r)$ for a unique $(\lambda,r)\in P^+\times\bn$. Moreover, $P(\lambda,r)$
is an indecomposable\index{indecomposable}  projective object of $\cal G$  and has $V(\lambda,r)$ as its unique irreducible quotient
 via a map $\tau_r\pi_\lambda$ which maps  $1\otimes v_\lambda\to v_\lambda$.\hfill\qedsymbol
\end{prop}
\noindent Since $\lie g=[\lie g,\lie g]$ one  deduces  that the kernel $K(\lambda,r)$ of the canonical surjection $P(\lambda,r)\twoheadrightarrow V(\lambda,r)$
is generated by $\tau_r\bu(\lie g[t])[1]\otimes_{\bu(\lie g)} V(\lambda)$, and if we note that $\bu(\lie g[t])[1]\cong(\lie g\otimes t)\bu(\lie g)$, we get
 $$K(\lambda,r)= \tau_r\left(\bu(\lie g[t])(\lie g\otimes t)\right)\otimes V(\lambda).$$
 A little more work proves that $$\Hom_{\cal G}(K(\lambda,r), V(\mu,s))\cong\begin{cases} 0,& s\ne r+1,\\ \Hom_{\lie g}(\lie g\otimes V(\lambda), V(\mu)),& s=r+1.\end{cases}$$
 We note the following
 \begin{cor} For  all $(\lambda,r),(\mu,s)\in P^+\times \bn$ we have
 \begin{align*}&\Ext^1_{\cal G}(V(\lambda,r), V(\mu,s))=0,
 \qquad s\ne r+1,\\ &\Ext^1_{\cal G}(V(\lambda,r), V(\mu,r+1))\cong  \Hom_{\lie g}(\lie g\otimes V(\lambda), V(\mu)).\tag*{\qed}
 \end{align*}
\end{cor}
 Motivated by the previous corollary, we define a partial order $\preceq$  on the index  set $P^+\times\bn$ of the simple modules in $\cal G$ by extending the cover relation $(\lambda,r)\prec(\mu,s)$ if and only if  $s=r+1$ and $\lambda-\mu\in\Phi\sqcup\{0\}$. Since a pair $(\mu,s)$ can cover only finitely many elements  it follows that  there exist only finitely many elements less than a given element. A subset $\Gamma$ of the poset $P^+\times\bn$ is called interval closed\index{interval closed} if for all $(\mu,s)\prec(\nu, k)\prec(\lambda,r)$ with $(\mu,s), (\lambda,r)\in\Gamma$ we have $(\nu,k)\in\Gamma$.

 \subsection{}\label{C:section 7.3} Even though an object $V$ of $\cal G$  could be infinite dimensional, one can still talk about the multiplicity of $V(\lambda,r)$
 in $V$ by setting $$[V:V(\lambda,r)]= [V_{\le r}: V(\lambda,r)], \qquad  V_{\le r}= \frac{V}{\bigoplus_{s>r}V[s]},$$
 where the right hand side of the preceding equality makes sense since $V_{\le r}$ is a finite--dimensional object of $\cal G$.  If $[V:V(\lambda,r)]>0$,
 then we call $V(\lambda,r)$ a composition factor\index{composition factor} of $V$.
 It follows from the definition of $P(\lambda,r)$ that $$[P(\lambda,r): V(\mu,s)]=\dim\Hom_{\lie g}(\bu(\lie g[t]_+)[s-r]\otimes V(\lambda), V(\mu)).$$
 Moreover, one can describe the $\lie g$-module structure of $\bu(\lie g[t]_+)$ in terms of symmetric powers of the adjoint representation of~$\lie g$.

 Given a finite interval closed\index{interval closed} subset $\Gamma$ of $P^+\times\bn$, let $\cal G[\Gamma]$ be the full subcategory of $\cal G$ consisting of modules $V$ whose composition factors\index{composition factor}   all lie in $\Gamma$, that is, $V\in\Ob\mathcal G[\Gamma]$ if and only if $$[V:V(\mu,s)]>0\implies (\mu,s)\in\Gamma.$$
Clearly, $\mathcal G[\Gamma]$ is closed under extensions, submodules and quotients.
 Given $V\in\Ob\cal G$, let $V_\Gamma$ be the maximal $\lie g[t]$--submodule of $V$  such that $$[V_\Gamma: V(\mu,s)]>0\implies (\mu,s)\notin\Gamma,$$ and set $V^\Gamma=V/V_\Gamma$. It is clear from the definition that $V^\Gamma\in\cal G[\Gamma]$ but it is generally not true that $V_\Gamma\in\cal G[\Gamma^c]$,
 where $\Gamma^c$ is the complement of $\Gamma$ in $P^+\times\bn$. However, as the following proposition shows, it is true for modules $P(\mu,s)$, with $(\mu,s)\in\Gamma$, and it is this fact that makes it possible to compute the ext quiver of the endomorphism algebra of the projective generator of the category $\cal G[\Gamma]$ for finite, interval closed\index{interval closed} $\Gamma$.
 \begin{prop}\label{C:cruxp}  Let $\Gamma$ be an interval closed\index{interval closed} subset of $\cal G$ and let $(\lambda,r), (\mu,s)\in\Gamma$. Then, $$[P(\lambda,r): V(\mu,s)]= [P(\lambda,r)^\Gamma: V(\mu,s)],$$ or equivalently $$ [P(\lambda,r)_\Gamma: V(\mu,s)] =0$$ and hence $P(\lambda,r)_\Gamma\in\cal G[\Gamma^c]$.\hfill\qedsymbol
 \end{prop}
\subsection{}\label{C:section 7.4} We continue to assume that $\Gamma$ is finite and interval closed\index{interval closed} and choose   $N\in\bn$  so that $(\mu,s)\in\Gamma$ only if $s<N$. Set \begin{gather*}P(\Gamma)=\bigoplus_{(\lambda,r)\in\Gamma}P(\lambda,r),\qquad \mathfrak A(\Gamma)=\End_{\cal G}(P(\Gamma)).  \end{gather*}
We have
$$\mathfrak A(\Gamma) \cong\End_{\cal G}(P(\Gamma)_{\le N})\cong\End_{\cal G[\Gamma]}(P(\Gamma)^\Gamma),$$ where the first isomorphism follows easily from the definitions and the second is an application of Proposition \ref{C:cruxp}. One then uses standard arguments to prove that
$\mathfrak A(\Gamma)$ is a finite--dimensional basic algebra whose left module category is equivalent to $\cal G[\Gamma]$.
Moreover, one can also prove that  $\mathfrak A(\Gamma)$  is a  directed quasi--hereditary algebra and an algorithm to compute the $\Ext$-quiver of this algebra along with the number of relations is given in \cite{CG1}. As $\Gamma$ varies over interval closed\index{interval closed} sets, one gets many interesting quivers and the algebras $\mathfrak A[\Gamma]$ have varying representation type.

We give one specific example, motivated by the
study of Kirillov-Reshetikhin modules, where the algebra is tame. Let $\lie g$ be of type $D_n$ with  $n\ge 6$ and use  the standard Bourbaki  numbering  of the Dynkin diagram. In particular, $4$ is not a spin node. Let $$\Gamma=\{(2\omega_4,0), (\omega_2+\omega_4, 1), (2\omega_2,2), (\omega_4, 2), (\omega_1+\omega_3,2), (\omega_2,3), (0,4)\}.$$ A simple exercise shows that $\Gamma$ is interval closed\index{interval closed}. With some further work, the results discussed so far can be used to prove that the category $\cal G[\Gamma]$  is equivalent to the module category of  the path algebra of the following quiver
 $$
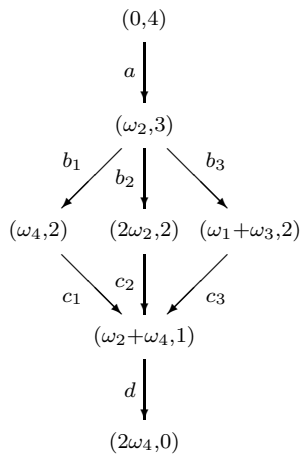
\makeatletter
\def\dggeometry{\unitlength=0.04pt\dg@YGRID=1\dg@XGRID=1}
\def\dgeverynode{\scriptstyle}
\begin{diagram}
\node{}\node{(0,4)}\arrow{s,l}{a}\\
\node{}\node{(\omega_2,3)}\arrow{se,l}{b_3}\arrow{s,l}{b_2}\arrow{sw,l}{b_1}\\
\node{(\omega_4,2)}\arrow{se,r}{c_1}\node{(2\omega_2,2)}\arrow{s,l}{c_2}\node{(\omega_1+\omega_3,2)}\arrow{sw,r}{c_3}\\
\node{}\node{(\omega_2+\omega_4,1)}\arrow{s,l}{d}\\
\node{}\node{(2\omega_4,0)}
\end{diagram}
$$
with  relations
$$
b_3a=0,\qquad dc_3=0,\qquad
c_1b_1+c_2b_2+c_3b_3=0.
$$
In particular, this algebra is quadratic. It can be shown to be
tame and of global dimension~$2$.

\subsection{}\label{C:section 7.5} We now further focus our study on finite--dimensional associative algebras which arise from the full subcategory $\cal G_2$ which consists of objects in $\cal G$ satisfying $$(\lie g\otimes t^r)V=0,\ \ r\ge 2.$$ One could define $\cal G_s$ for $s\in\bn$ in the obvious way. For $s=1$ it is easy (and has been said elsewhere in different forms throughout these notes) to see that the category $\cal G_1$ is semi-simple. The irreducible objects are just $\ev^*_0V(\lambda)$ for $\lambda\in P^+$. Thus $\cal G_2$ is the first interesting case which is not semisimple, making it reasonable to limit our attention to this category.
It is worth mentioning that restricted Kirillov-Reshetikhin modules in classical types are objects in this category.
An object in $\cal G_2$ is actually a representation of the graded truncated  Lie algebra $\lie g\otimes(\bc[t]/(t^2))$.
This algebra is  isomorphic to the semi--direct product $\lie g\ltimes\lie g_{\ad}$ where $\lie g_{\ad}$ is the adjoint representation of $\lie g$. Although this algebra is supported only in grades zero and one, the universal enveloping algebra is $\bn$--graded and we have an isomorphism of associative algebras $$\bu(\lie g\ltimes\lie g_{\ad})\cong S(\lie g_{\ad})\otimes \bu(\lie g),$$ and the $r^{th}$ graded piece is $$\bu(\lie g\ltimes\lie g_{\ad})[r]\cong S^r(\lie g_{\ad})\otimes \bu(\lie g).$$ Here $S(\lie g_{\ad})$ is the symmetric algebra of $\lie g_{\ad}$ and $S^r(\lie g_{\ad})$ is the $r^{th}$--symmetric space. One can now use the Koszul complex of $S(\lie g)$ to construct an explicit projective resolution for any simple object in~$\mathcal G_2$ and to compute all $\Ext$ spaces
between simple objects in $\cal G_2$.
\begin{prop} \label{C:extj} For all $(\lambda,r),(\mu,s)\in P^+\times\bn$, we have
\begin{equation*}
\Ext^j_{\cal G_2}(V(\lambda,r), V(\mu,s))\cong\begin{cases} 0,&j\ne s-r,\\ \Hom_{\lie g}(\bigwedge^j\lie g_{\ad}\otimes V(\lambda), V(\mu)),&j=s-r.\end{cases}
\end{equation*}
In particular, $\Ext^j_{\cal G_2}(M,N)=0$ if $M,N$ are finite dimensional and $j>\dim\lie g$.\qed
\end{prop}
This proposition strongly suggests that there is an associative algebra with a Koszul grading in the background and we now describe how to find this algebra.

\subsection{}\label{C:section 7.6} We first need to refine the partial ordering $\preceq$ on $P^+\times\bn$. Instead of allowing all possible roots in the cover relation, we restrict ourselves to particular subsets of $\Phi^+$ obtained as follows. Given $\psi\in P$, let
$$
\Phi^+_\psi=\{\alpha\in R^+:\kappa(\psi,\alpha)\ge \kappa(\psi, \beta)\,\text{for all $\beta\in \Phi$}\},
$$
and define a cover relation  by $(\mu,s)$ covers $(\lambda,r)$ iff $s=r+1$ and $\lambda-\mu\in \Phi^+_\psi$. Let $\le_\psi$ be the corresponding partial order or, equivalently, $(\lambda, r)\le_\psi(\mu,s)$ iff $$\mu-\lambda=\sum_{\alpha\in \Phi^+_\psi} n_\alpha\alpha,\quad n_\alpha\in\bn,   \qquad
\sum_{\alpha\in \Phi^+_\psi}n_\alpha =s-r.$$ The definition of $\Phi^+_\psi$ ensures that this order is well--defined. The following was proved in \cite{CG2}:
\begin{thm} Let $(\lambda,r)\in P^+\times \bn$ and let $\Gamma=\{(\mu,s): (\mu,s)\le_\psi(\lambda,r)$. Then the endomorphism algebra $\mathfrak A[\Gamma]$ of the projective generator of $\cal G_2[\Gamma]$ admits a natural Koszul grading and
is of global dimension at most the cardinality of $\Phi^+_\psi$. Moreover the maximal dimension is attained for a suitable choice of $(\lambda,r)$.\qed
\end{thm}
One can define certain limits of the algebras as $\Gamma$ varies and construct infinite--dimensional Koszul algebras $\mathfrak A_\psi$ of left global dimension
$|\Phi^+_\psi|$. The Koszul dual of these algebras are also studied in \cite{CG2}, although one does not understand the module category of the dual in the context of representations of Lie algebras. We conclude these notes with examples of the quivers associated to the infinite--dimensional Koszul algebras arising from our study. Further details may be found in \cite{G}.

 \subsection{}\label{C:section 7.7} Suppose that~$\lie g$ is not of type~$A$ or~$C$ so that
 there exists a unique~$i_0\in I$ such that~$\theta-\alpha_{i_0}\in \Phi^+$ and then,
$$\Phi^+_{\omega_{i_0}}=\{\theta,\theta-\alpha_{i_0}\}.$$  Then every connected subalgebra of~$\mathfrak A_{\omega_{i_0}}$ is isomorphic to the path algebra of the translation quiver
\begin{equation}\label{SP100.10}
\def\dgeverynode{\scriptscriptstyle}\divide\dgARROWLENGTH by 2\dgVERTPAD=2pt\dgHORIZPAD=2pt
\begin{diagram}
\node[3]{}\node{\vdots}\arrow{s}\node{\vdots}\arrow{s}
\\
\node{}\node{}\node{(0,2)}\arrow{s}\node{(1,2)}\arrow{w}\arrow{s}\node{(2,2)}\arrow{w}\arrow{s}\node{\displaystyle\cdots}\arrow{w}
\\
\node{}\node{(0,1)}\arrow{s}\node{(1,1)}\arrow{w}\arrow{s}\node{(2,1)}\arrow{w}\arrow{s}\node{(3,1)}\arrow{w}\arrow{s}\node{\displaystyle\cdots}\arrow{w}
\\
\node{(0,0)}\node{(1,0)}\arrow{w}\node{(2,0)}\arrow{w}\node{(3,0)}\arrow{w}\node{(4,0)}\arrow{w}\node{\displaystyle\cdots}\arrow{w}
\end{diagram}
\end{equation}
with the mesh relations as indicated.

\subsection{}\label{C:section 7.8} Let  $\lie g$ be of type~$C$. Then $\Phi^+_{\omega_2}=\{\theta,\theta-\alpha_1,\theta-2\alpha_1\}$ and the quiver of $\mathfrak A_{\omega_2}$ is
the (infinite, if rank of $\lie g$ is greater than $2$) union of connected components of the following two types:
$$
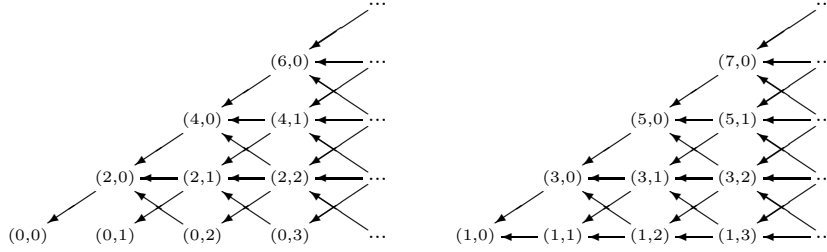
\def\dgeverynode{\scriptscriptstyle}
\divide\dgARROWLENGTH by 2\dgVERTPAD=4pt\dgHORIZPAD=4pt
\begin{diagram}
\node{}\node{}\node{}\node{}\node{\scriptstyle\dots}\arrow{sw}
\\
\node{}\node{}\node{}\node{(6,0)}\arrow{sw}\node{\scriptstyle\dots}\arrow{w}\arrow{sw}
\\
\node{}\node{}\node{(4,0)}\arrow{sw}\node{(4,1)}\arrow{w}\arrow{sw}\node{\scriptstyle\dots}\arrow{nw}\arrow{sw}\arrow{w}
\\
\node{}\node{(2,0)}\arrow{sw}\node{(2,1)}\arrow{w}\arrow{sw}\node{(2,2)}\arrow{nw}\arrow{w}\arrow{sw}
\node{\scriptstyle\dots}\arrow{nw}\arrow{sw}\arrow{w}\\
\node{(0,0)}\node{(0,1)}\node{(0,2)}\arrow{nw}\node{(0,3)}\arrow{nw}\node{\scriptstyle\dots}\arrow{nw}\\
\end{diagram}
\qquad
\begin{diagram}
\node{}\node{}\node{}\node{}\node{\scriptstyle\dots}\arrow{sw}
\\
\node{}\node{}\node{}\node{(7,0)}\arrow{sw}\node{\scriptstyle\dots}\arrow{w}\arrow{sw}
\\
\node{}\node{}\node{(5,0)}\arrow{sw}\node{(5,1)}\arrow{w}\arrow{sw}\node{\scriptstyle\dots}\arrow{nw}\arrow{sw}\arrow{w}
\\
\node{}\node{(3,0)}\arrow{sw}\node{(3,1)}\arrow{w}\arrow{sw}\node{(3,2)}\arrow{nw}\arrow{w}\arrow{sw}
\node{\scriptstyle\dots}\arrow{nw}\arrow{sw}\arrow{w}\\
\node{(1,0)}\node{(1,1)}\arrow{w}\node{(1,2)}\arrow{nw}\arrow{w}\node{(1,3)}\arrow{nw}\arrow{w}\node{\scriptstyle\dots}\arrow{nw}\arrow{w}\\
\end{diagram}
$$
Both are translation quivers with~$\tau((m,n))=(m,n-2)$, $m>0$, $n\ge 2$. The relations are: the commutativity relations in
$$
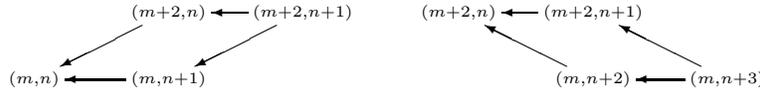
\def\dgeverynode{\scriptscriptstyle}
\divide\dgARROWLENGTH by 2
\dgVERTPAD=4pt\dgHORIZPAD=4pt
\begin{diagram}
\node{}\node{(m+2,n)}\arrow{sw}\node{(m+2,n+1)}\arrow{w}\arrow{sw}\\
\node{(m,n)}\node{(m,n+1)}\arrow{w}
\end{diagram}\qquad
\begin{diagram}
\node{(m+2,n)}\node{(m+2,n+1)}\arrow{w}\\\node{}\node{(m,n+2)}\arrow{nw}\node{(m,n+3)}\arrow{w}\arrow{nw}
\end{diagram}
$$
for all~$m>0$, $n\in\bz_+$, the zero relations $(2,n)\from (2,n+1)\from (0,n+3)$, $n\ge 0$
and
\begin{multline*}
m^2((m,n)\from (m+2,n)\from (m,n+2))-(m+2)^2((m,n)\from (m-2,n+2)\from (m,n+2))\\+(m+1)((m,n)\from(m,n+1)\from(m,n+2)),\qquad m>1,
\end{multline*}
and, finally, $((1,n)\from (3,n)\from (1,n+2))+2((1,n)\from(1,n+1)\from(1,n+2))$.
Thus, if~$\ell=2$ and~$|\Phi^+_\psi|>1$, the algebra $\mathfrak A_\psi$ is the direct sum of two non-isomorphic connected Koszul subalgebras of left global dimension~$3$.

\subsection{}\label{C:section 7.9} Let $\lie g$ be of type~$G_2$ and suppose that $\alpha_1$ is the long simple root. Then $\Phi^+_{\omega_1-\omega_2}=\{\alpha_1,\theta\}$
and the algebra
$\mathfrak A_{\omega_1-\omega_2}$ is the direct sum of three isomorphic connected subalgebras with quivers
$$
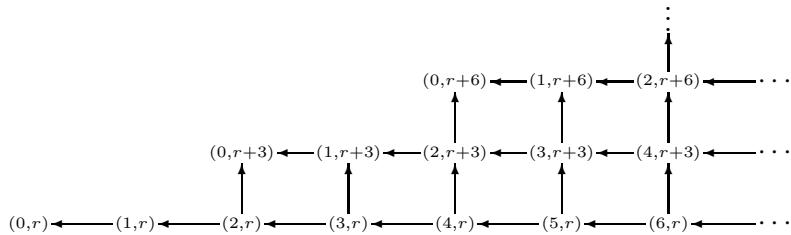
\def\dgeverynode{\scriptscriptstyle}
\divide\dgARROWLENGTH by 2\dgVERTPAD=2pt\dgHORIZPAD=2pt
\begin{diagram}
\node[6]{}\node{\vdots}
\\
\node[4]{}\node{(0,r+6)}\node{(1,r+6)}\arrow{w}\node{(2,r+6)}\arrow{n}\arrow{w}\node{\displaystyle\cdots}\arrow{w} \\
\node{}\node{}\node{(0,r+3)}\node{(1,r+3)}\arrow{w}\node{(2,r+3)}\arrow{n}\arrow{w}\node{(3,r+3)}\arrow{w}\arrow{n}\node{(4,r+3)}\arrow{w}\arrow{n}
\node{\displaystyle\cdots}\arrow{w}\\
\node{(0,r)}\node{(1,r)}\arrow{w}\node{(2,r)}\arrow{n}\arrow{w}\node{(3,r)}\arrow{w}\arrow{n}\node{(4,r)}\arrow{w}\arrow{n}
\node{(5,r)}\arrow{w}\arrow{n}\node{(6,r)}\arrow{w}\arrow{n}\node{\displaystyle\cdots}\arrow{w}
\end{diagram}
$$
where $0\le r\le 2$. The relations are again the mesh relations, the translation map being $$\tau((m,3k+r))=(m+3,3(k-1)+r),\ \ m\in\bz_+,\ \  k>0.$$

\subsection{}\label{C:section 7.10} In conclusion, we note that there is some work  \cite{BBS}, \cite{Bil}, \cite{CLE}, \cite{La}, \cite{Moo-Rao-Yok} on the representation theory of extended affine algebras where the full center does not act trivially. There is also the recent work \cite{FZ} in which the authors give a representation theoretic way to interpret the third cohomology classes of the  double loop algebras. It is perhaps clear by now, that there are many interesting avenues to pursue in the representation theory of affine, toroidal and extended affine Lie algebras and while the references are by no means anywhere near exhaustive they should provide the reader with some pointers to the literature.

\medskip

{\small
\begin{flushright}
{\bf Vyjayanthi Chari}\\
 Department of Mathematics\\
 University of California\\
Riverside, CA 92521 USA\\
{\tt chari@math.ucr.edu}
\end{flushright}}

\vspace{4mm} \noindent \small{Vyjayanthi Chari: {\it Department of
Mathematics, University of California at Riverside, CA 92521, USA.}}
\end{document}